\everymath{\displaystyle}
\documentclass{article}

\usepackage[english]{babel}
\usepackage[T1]{fontenc}

\usepackage{amsmath}
\usepackage{amssymb}
\usepackage{amsthm}
\usepackage{mathtools}
\usepackage{physics}
\usepackage{array}
\usepackage[per-mode=symbol,
            inter-unit-product=\cdot,
            separate-uncertainty=true]{siunitx}
\usepackage[margin=1cm]{caption}
\usepackage{graphicx}
\usepackage{titlesec}
\usepackage[hypertexnames=false,
            colorlinks=true,
            pdfborderstyle={/S/U/W 0},
            citecolor=blue]{hyperref}
\usepackage[capitalise,noabbrev]{cleveref}
\usepackage{bookmark}
\usepackage[a4paper, total={6.5in, 8.66in}]{geometry}
\usepackage{multirow}
\usepackage{todonotes}
\usepackage{subcaption}
\usepackage[toc]{appendix}
\usepackage{cases}
\usepackage{multicol}
\usepackage{multirow}
\usepackage{booktabs}

\newcommand{\fluid}{_\mathrm f} 
\newcommand{\solid}{_\mathrm s} 
\newcommand{\ale}{_\mathrm{ALE}} 

\newcommand{\lv}{^\mathrm{LV}} 
\newcommand{\la}{^\mathrm{LA}} 
\newcommand{\ao}{^\mathrm{AA}} 
\newcommand{\ring}{^\mathrm{ring}} 

\newcommand{\ar}{_\mathrm{AR}}
\newcommand{\ven}{_\mathrm{VEN}}
\newcommand{\sys}{^\mathrm{SYS}}
\newcommand{\pul}{^\mathrm{PUL}}

\newcommand{\lifex}{\texttt{life\textsuperscript{x}} }
\newcommand{\dealii}{\texttt{deal.II} }

\renewcommand{\phi}{\varphi} 
\renewcommand{\epsilon}{\varepsilon} 
\renewcommand{\laplacian}{\Delta} 
\renewcommand{\hat}{\widehat} 

\DeclareSIUnit\mmhg{mmHg}

\newlength{\figwidth}
\title{A mathematical model that integrates cardiac electrophysiology, mechanics and fluid dynamics: application to the human left heart}
\author{Michele Bucelli$^{1, *}$,
        Alberto Zingaro$^1$,
        Pasquale Claudio Africa$^1$, \\
        Ivan Fumagalli$^1$,
        Luca Dede'$^1$,
        Alfio Quarteroni$^{1, 2}$}

\date{\footnotesize \textsuperscript{1} MOX, Department of Mathematics, Politecnico di Milano, P.zza Leonardo da Vinci 32, 20133 Milan, Italy \\
      \textsuperscript{2} Mathematics Institute, EPFL, Av. Piccard, CH-1015 Lausanne, Switzerland (Professor Emeritus) \\[2ex]%
      \today}

\begin{document}

\maketitle
\makeatletter{\renewcommand*{\@makefnmark}{}\footnotetext{* Corresponding author, email address: \texttt{michele.bucelli@polimi.it}}\makeatother}

\begin{abstract}
    We propose a mathematical and numerical model for the simulation of the heart function that couples cardiac electrophysiology, active and passive mechanics and hemodynamics, and includes reduced models for cardiac valves and the circulatory system. Our model accounts for the major feedback effects among the different processes that characterize the heart function, including electro-mechanical and mechano-electrical feedback as well as force-strain and force-velocity relationships. Moreover, it provides a three-dimensional representation of both the cardiac muscle and the hemodynamics, coupled in a fluid-structure interaction (FSI) model. By leveraging the multiphysics nature of the problem, we discretize it in time with a segregated electrophysiology-force generation-FSI approach, allowing for efficiency and flexibility in the numerical solution. We employ a monolithic approach for the numerical discretization of the FSI problem. We use finite elements for the spatial discretization of those partial differential equations that contribute to the model. We carry out a numerical simulation on a realistic human left heart model, obtaining results that are qualitatively and quantitatively in agreement with physiological ranges and medical images.
\end{abstract}

{\textbf{Keywords:} multiphysics modeling, cardiac modeling, electromechanics, fluid-structure interaction, blood circulation}

\section{Introduction}

\begin{figure}
    \centering
    \begin{subfigure}{0.4\textwidth}
        \centering
        \includegraphics[width=\textwidth]{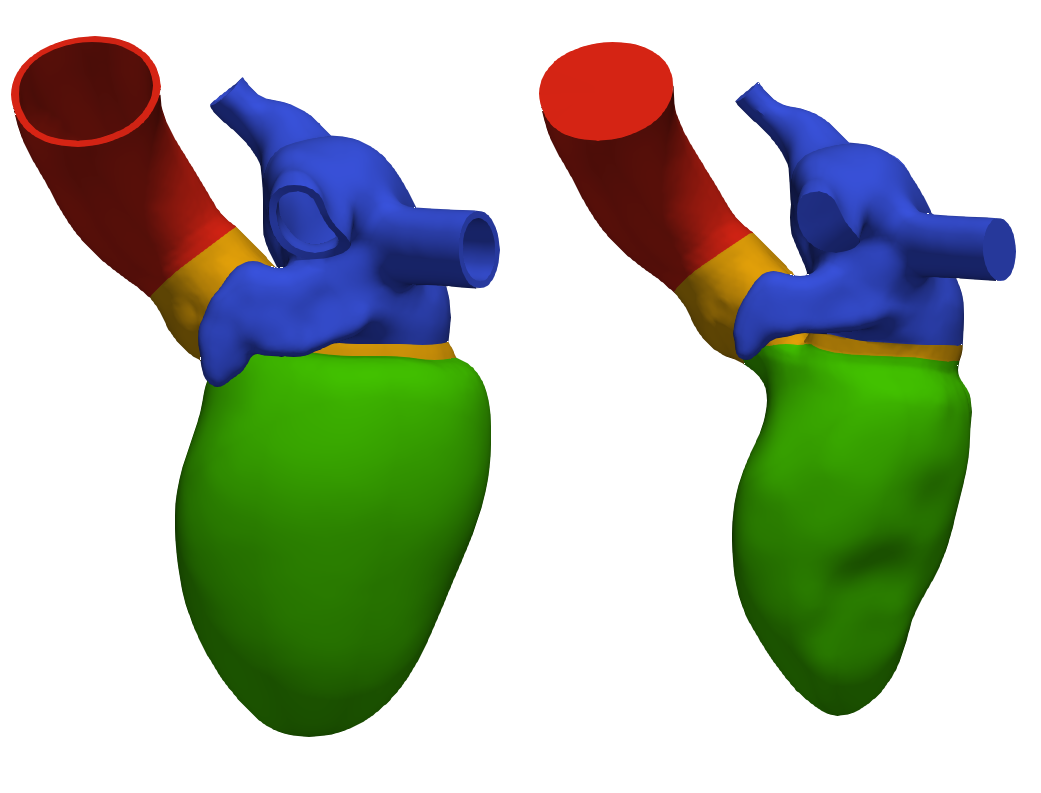} \\
        \includegraphics[width=0.35\textwidth]{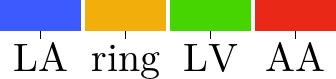}
        \caption{}
    \end{subfigure}
    \begin{subfigure}{0.4\textwidth}
        \centering
        \includegraphics[width=\textwidth]{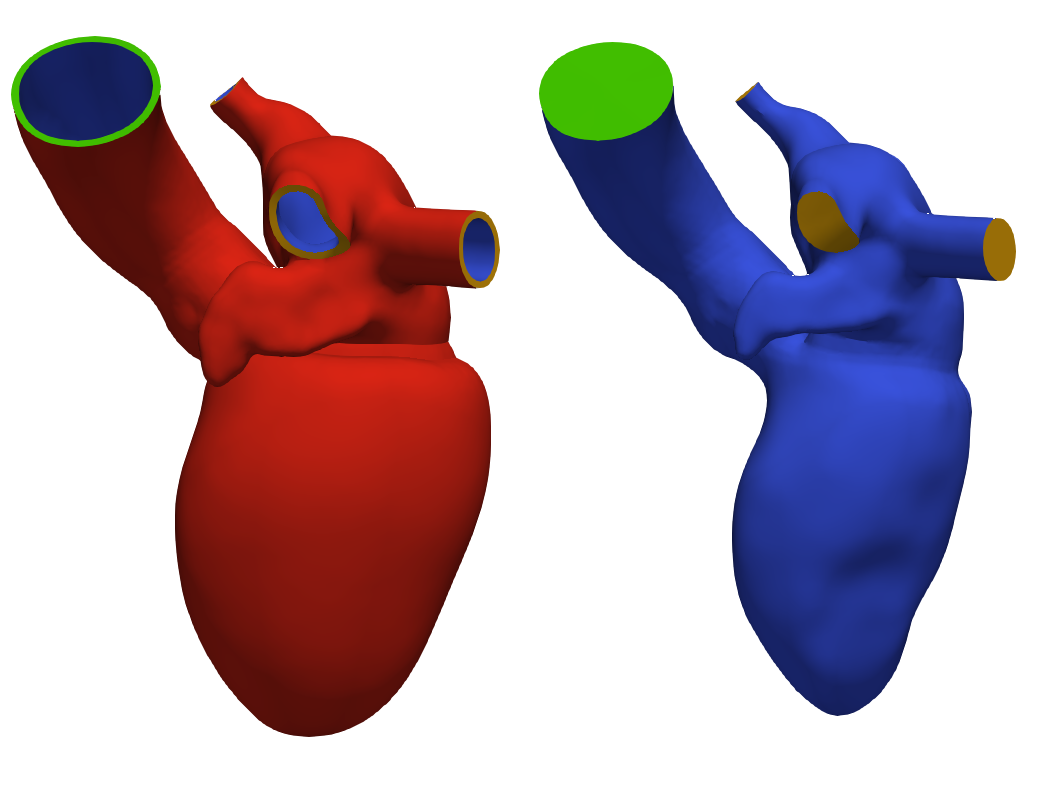} \\
        \includegraphics[width=0.35\textwidth]{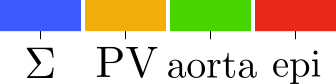}
        \caption{}
        \label{fig:domain-boundary}
    \end{subfigure}

    \caption{Solid (left) and fluid (right) computational domains. Colors represent the partition into $\Omega_\text{i}\lv$, $\Omega_\text{i}\la$, $\Omega_\text{i}\ao$ and $\Omega_\text{i}\ring$ in (a), with $\text{i} \in \{\text{s}, \text{f}\}$, and the different boundary portions in (b).}
    \label{fig:domain}
\end{figure}

Cardiovascular diseases represent the major cause of death in the adult population of the western world \cite{benjamin2017heart}. While in-vivo measuring techniques allow to inspect and quantify the heart function and dysfunction, these measures often lack resolution or accuracy, and may be invasive. Mathematical models for computational medicine can provide tools to simulate the human heart function, complementing experimental measurements, providing further insight into the cardiac function, and assisting in the development of personalized treatments \cite{gerach2021electro,gray2018patient,karabelas2022global,niederer2019computational,quarteroni2017integrated,quarteroni2019mathematical}.

The heart function results from the interplay of several different physical processes \cite{katz2010physiology,klabunde2011cardiovascular}, ranging from the sub-cellular scale to the tissue one. Electrochemical processes drive the excitation of cardiac muscular cells, resulting in the generation of an active force during contraction that, together with the passive mechanical properties of the cardiac tissue, interacts with the blood within the heart chambers, pumping it towards the circulatory system. Understanding and then modeling the multiphysics and multiscale nature of the heart function, as well as feedback effects among its components, is crucial for constructing an accurate computational model of the heart \cite{hosoi2010multi,okada2019clinical,quarteroni2017integrated,santiago2018fully}.

Several mathematical and numerical heart models have been proposed in recent years. Most of them focus on some specific feature of the heart function, by surrogating the remaining ones with models of reduced dimensionality:
electrophysiology \cite{arevalo2016arrhythmia,bucelli2021multipatch,del2022fast,gillette2021framework,piersanti2021modeling,romero2010effects,trayanova2011whole,vergara2016coupled},
cardiac mechanics and electromechanics \cite{augustin2016patient,baillargeon2014living,dede2020segregated,fedele2022comprehensive,gerach2021electro,gerbi2019monolithic,gurev2011models,levrero2020sensitivity,pfaller2019importance,piersanti20223d,regazzoni2022cardiac,salvador2021electromechanical,strocchi2020simulating,usyk2002computational}
or computational fluid dynamics (CFD) of the blood \cite{chnafa2014image,collia2019analysis,karabelas2022global,kronborg2022computational,rigatelli2021applications,spuhler2020high,terahara2022computational,this2020augmented,vedula2016effect,zingaro2021hemodynamics,zingaro2022geometric}.
Only few works also consider the interplay between hemodynamics and cardiac mechanics in a fluid-structure interaction (FSI) framework \cite{cheng2005fluid,einstein2010fluid,feng2019analysis,gao2017coupled,hirschhorn2020fluid}, while neglecting or simplifying the electrical processes generating the contraction \cite{brenneisen2021sequential,khodaei2021personalized,nordsletten2011fluid,zhang2001analysis}.

Albeit each of the above mentioned models can provide meaningful insight into the cardiac function in both healthy and pathological conditions, they often neglect the feedback mechanisms that relate the different components. On the other hand, fully integrated models featuring multiphysics coupling of electrophysiology, mechanics and fluid dynamics \cite{bucelli2022stable,gerbi2018numerical,santiago2018fully,sugiura2022ut,verzicco2022electro} can provide a very accurate description of the physics of the heart, at the price of a high model complexity and large computational cost.

One of the first fully coupled heart models was provided by the UT-Heart simulator \cite{hosoi2010multi,okada2019clinical,sugiura2022ut}, and was subsequently employed for personalized clinical case studies \cite{kariya2020personalized}. Another fully coupled model was proposed in \cite{santiago2018fully}. There, the authors focus on right and left ventricular systole, while providing a simplified description of the atria and neglecting the presence of valves. In \cite{viola2020fluid,viola2021effects}, a coupled fluid-structure-electrophysiology left heart model is presented, relying on the immersed boundary method in a finite elements/finite differences combined framework for the numerical discretization. The model uses a simplified phenomenological description of force generation \cite{nash2004electromechanical}, neglecting feedbacks from fiber shortening; moreover, the interplay between heart and circulation is treated in a simplified approach. The same authors also provide a GPU-accelerated version of the same model in \cite{viola2022fsei} for computational speedup. Finally, in \cite{battista2017bifurcations}, a simplified 2D fluid-structure-electrophysiology interaction model is developed for embryonic hearts. To the best of our knowledge, no other works focus on a fully integrated model of the heart.

In this work, we propose a novel multiphysics coupled model featuring a three-dimensional description of cardiac electrophysiology, active and passive mechanics, and hemodynamics, together with a closed-loop lumped-parameter model to simulate systemic and pulmonary circulation. We will refer to the different physical components (electrophysiology, contractile force generation, mechanics, fluid dynamics, circulation) as \textit{core models}. The proposed model includes several of the feedback mechanisms that regulate the heart function, including mechano-electric feedbacks \cite{salvador2022role}, force-strain and force-velocity relationships \cite{regazzoni2020biophysically}, FSI between the blood and the  muscle \cite{bucelli2022partitioned}, and the feedback between the heart and the circulation \cite{hirschvogel2017monolithic,regazzoni2022cardiac,zingaro2022geometric}. In particular, the feedback between force generation and fiber shortening and shortening velocity was found to be highly important in cardiac electromechanics \cite{fedele2022comprehensive}. Moreover, FSI modeling allows to capture dynamic effects such as the presence of pressure gradients and pressure waves, overcoming some of the limitations of uncoupled CFD and electromechanics models \cite{this2020augmented,zingaro2022modeling}.
We test our model on a realistic human left heart, comprising the left ventricle, left atrium and a portion of the ascending aorta, as well as the mitral and aortic valves.

Due to the complexity and large scale of the problem at hand, it is crucial to employ efficient numerical schemes for its solution. A critical issue is the numerical treatment of the coupling conditions between the different core models. One possible approach is to solve the fully coupled problem monolithically \cite{gerbi2018numerical}. However, this requires the development of complex solvers and ad-hoc preconditioners, and it lacks flexibility in the choice of the discretization schemes and parameters for each core model. Inspired by \cite{regazzoni2022cardiac}, we choose a segregated-staggered scheme, in which all coupling terms are treated explicitly when possible. Implicit, monolithic coupling is used for fluid-structure interaction \cite{bucelli2022partitioned}. Spatial discretization is achieved by means of the finite element method \cite{hughes2012finite,quarteroni2017numerical}. The proposed computational framework leverages high-performance computing techniques to enable large-scale simulations, based on the finite element library \lifex \cite{africa2022lifexcore,africa2022lifex,lifex}.

We run a numerical simulation and compare the value of several indicators against normal ranges, obtaining a satisfactory agreement. A qualitative match with physiological behavior is also observed in terms of ventricular pressure and volume traces, as well as three dimensional deformation and flow patterns. Overall, the results indicate that the proposed computational framework can provide an accurate and physiologically sound description of the physics of the heart.

This paper is structured as follows: in \cref{sec:mathematical}, we present in detail the mathematical models employed. \Cref{sec:numerical} describes the numerical scheme, and \cref{sec:results} presents and analyzes the simulation results. Finally, in \cref{sec:limitations}, we point at some shortcomings in the proposed model, and, in \cref{sec:conclusions}, we draw some conclusive remarks.

\section{Mathematical models}
\label{sec:mathematical}

Let $\Omega \subset \mathbb{R}^3$ be an open, bounded domain, changing in time and representing the volume occupied by a human left heart at every time instant during the heartbeat. We partition the domain into the subdomains $\Omega\fluid$ and $\Omega\solid$ (see \cref{fig:domain}), representing the volume occupied by the blood and that occupied by the solid heart wall respectively. We denote their interface by $\Sigma = \partial\Omega\fluid\cap\partial\Omega\solid$. We further partition the fluid and solid subdomains into $\Omega_\mathrm{i}\lv$, $\Omega_\mathrm{i}\la$, $\Omega_\mathrm{i}\ao$ and $\Omega_\mathrm{i}\ring$, for $\mathrm{i} \in \{\mathrm{f}, \mathrm{s}\}$, representing the left ventricle (LV), the left atrium (LA), the ascending aorta (AA) and the atrioventricular ring, for either the fluid or the solid domains. On the boundary of any of the defined sets, $\mathbf n$ denotes the outward directed normal unit vector. On $\Sigma$, $\mathbf n$ denotes the normal unit vector directed outward from the fluid domain and inward into the solid domain.

All the previously defined domains are moving in time as the heart beats (to keep the notation light, we omit the explicit dependence of the domains on time). We define a fixed \textit{reference configuration} $\hat\Omega$ to keep track of the deformation. For each of the previously defined sets, we use a hat to denote the corresponding set in reference configuration. We refer to the moving configuration as the \textit{current configuration}. The evolution in time of the current configuration is described by the following maps:
\begin{align*}
    \mathcal{L}\solid &: \hat\Omega\solid \cross(0,T) \to \Omega\solid \qquad\qquad \Omega\solid = \{\mathbf x = \mathcal{L}\solid(\hat{\mathbf x}, t)\;, \hat{\mathbf x}\in\hat\Omega\solid\}\;, \\
    \mathcal{L}\fluid &: \hat\Omega\fluid \cross(0,T) \to \Omega\fluid \qquad\qquad \Omega\fluid = \{\mathbf x = \mathcal{L}\fluid(\hat{\mathbf x}, t)\;, \hat{\mathbf x}\in\hat\Omega\fluid\}\;.
\end{align*}
The precise definition of the maps depends on the physical models defined on each subdomain, and is detailed in the following sections.

We denote by $t \in (0, T)$ the independent variable representing time, by $\mathbf x$ the spatial coordinates in the current configuration, and by $\hat{\mathbf x}$ the spatial coordinates in reference configuration.

We consider a coupled problem involving several physical models: electrophysiology, active force generation, cardiac FSI and circulation hemodynamics. The unknowns of the model are determined by solving a coupled system of differential equations representing an \textit{electrophysiology-mechanics-fluid dynamics} (EMF) interaction problem:
\begin{align*}
    & v: \hat\Omega\solid\lv \times (0, T) \to \mathbb{R} & \text{transmembrane potential,} \\
    & \mathbf w: \hat\Omega\solid\lv \times (0, T) \to \mathbb{R}^{N_\mathrm{ion}^{\mathbf w}} & \text{ionic gating variables,} \\
    & \mathbf z: \hat\Omega\solid\lv \times (0, T) \to \mathbb{R}^{N_\mathrm{ion}^{\mathbf z}} & \text{ionic concentrations,} \\
    & \mathbf s: \hat\Omega\solid\lv \times (0, T) \to \mathbb{R}^{N_\mathrm{act}} & \text{activation state,} \\
    & \mathbf d: \hat\Omega\solid \times (0, T) \to \mathbb{R}^3 & \text{solid displacement,} \\
    & \mathbf d\ale: \hat\Omega\fluid \times (0, T) \to \mathbb{R}^3 & \text{fluid domain displacement,} \\
    & \mathbf u: \Omega\fluid \times (0, T) \to \mathbb{R}^3 & \text{fluid velocity,} \\
    & p: \Omega\fluid \times (0, T) \to \mathbb{R} & \text{pressure,} \\
    & \mathbf c: (0, T) \to \mathbb{R}^{N_\mathrm{circ}} & \text{circulation state variables.}
\end{align*}

The remainder of this section is devoted to the description of the models and equations that compose the EMF problem. In the following, we omit the numerical values of all the parameters mentioned, that can be found in \cref{app:parameters}.

\subsection{Fiber generation}

Cardiac tissue is organized in sheets of fibers, which determine a preferential direction for the conduction of the electrical signal \cite{katz2010physiology,klabunde2011cardiovascular,roberts1979influence}, and are relevant in characterizing the passive and active mechanical properties of the cardiac muscle \cite{eriksson2013influence,gil2019influence}. We incorporate fibers in the model by defining, at every point $\hat{\mathbf x} \in \hat\Omega\solid$, a local orthonormal basis $\{\mathbf f_0, \mathbf s_0, \mathbf n_0\}$ representing the local direction of fibers, fiber sheets and normal to fiber sheets \cite{piersanti2021modeling}.

As the fiber direction is rarely available from experimental data, it is often reconstructed in a preprocessing step using rule-based algorithms \cite{doste2019rule,piersanti2021modeling,rossi2014thermodynamically}. For our left heart model, we combine different Laplace-Dirichlet Rule-Based Methods (LDRBMs), namely the one presented in \cite{bayer2012novel} for ventricular fibers and the one presented in \cite{piersanti2021modeling} for left atrial fibers.

\subsection{Electrophysiology}
\label{sec:electrophysiology}

Cardiac cells can be excited by an electrical stimulus, triggering a series of subcellular mechanisms resulting in an \textit{action potential} that manifests as a variation in time of the \textit{transmembrane potential} \cite{collifranzone2014mathematical,katz2010physiology,klabunde2011cardiovascular,sundnes2007computing}, i.e. the difference of potential between the intracellular and extracellular spaces.

We model the ventricle $\hat\Omega\solid\lv$ as electrically excitable, while other portions of the solid domain, including $\hat\Omega\solid\la$, are regarded as electrically passive, i.e. they do not generate action potentials. We remark that this is a simplification for what concerns the LA (see \cref{sec:limitations}).

The evolution of the transmembrane potential $v$ is described by the monodomain equation \cite{collifranzone2014mathematical} augmented with the mechano-electric feedbacks \cite{salvador2022role}. Let $F = I + \grad\mathbf d$ and $J = \det F$. The monodomain equation reads
\begin{equation}
    \begin{dcases}
        J \chi C_\text{m} \pdv{v}{t} - \div(J F^{-1} D_\mathrm{m} F^{-T} \grad v) + J \chi I_\mathrm{ion}(v, \mathbf w, \mathbf z) = J \chi I_\mathrm{app}(\hat{\mathbf x}, t) & \text{in }\hat\Omega\solid\lv \times (0, T)\;, \\
        J F^{-1} D_\mathrm{m} F^{-T}\grad v \cdot \mathbf n = 0 & \text{on }\partial\hat\Omega\solid\lv \times (0, T)\;, \\
        v = v_0 & \text{in }\hat\Omega\solid\lv \times \{0\}\;,
    \end{dcases}
    \label{eq:monodomain}
\end{equation}
where $\chi$ and $C_\text{m}$ are the membrane surface-to-volume ratio and membrane capacitance, respectively.

The vector $\mathbf w$ collects the recovery (or gating) variables \cite{collifranzone2014mathematical}, that for a single cell represent fractions of open ionic channels, while $\mathbf z$ is a vector of ionic concentrations. Most notably, one of the variables within the vector $\mathbf z$ represents the intracellular calcium concentration $[\mathrm{Ca}^{2+}]_\mathrm{i}$. The evolution of $\mathbf w$ and $\mathbf z$ is modeled by coupling \eqref{eq:monodomain} with the ionic model by Ten Tusscher and Panfilov \cite{ten2006alternans}, which also defines the ionic current $I_\mathrm{ion}(v, \mathbf w, \mathbf z)$. The model is expressed by a system of ODEs defined at each point $\hat{\mathbf x} \in \hat\Omega\solid\lv$:
\begin{equation*}
    \begin{dcases}
        \pdv{\mathbf w}{t} = \mathbf F_\mathrm{ion}^\mathbf{w}(v, \mathbf w) & \text{in } \hat\Omega\solid\lv \times (0, T)\;, \\
        \pdv{\mathbf z}{t} = \mathbf F_\mathrm{ion}^\mathbf{z}(v, \mathbf w, \mathbf z) & \text{in } \hat\Omega\solid\lv \times (0, T)\;, \\
        \mathbf w = \mathbf w_0 & \text{in } \hat\Omega\solid\lv \times \{0\}\;, \\
        \mathbf z = \mathbf z_0 & \text{in } \hat\Omega\solid\lv \times \{0\}\;.
    \end{dcases}
\end{equation*}
We remark that $\mathbf F_\mathrm{ion}^\mathbf{w}$ and $\mathbf F_\mathrm{ion}^\mathbf{z}$ do not involve spatial derivatives of $\mathbf w$, $\mathbf z$ or $v$. For the precise definition of $\mathbf F_\mathrm{ion}^\mathbf{w}$, $\mathbf F_\mathrm{ion}^\mathbf{z}$ and $I_\mathrm{ion}$, we refer to \cite{ten2006alternans}. The initial states $v_0$, $\mathbf w_0$ and $\mathbf z_0$ are obtained by solving a reduced, zero-dimensional monodomain equation for a large number of heartbeats, until a periodic limit cycle is reached \cite{collifranzone2014mathematical,regazzoni2022cardiac}.

The tensor $D_\mathrm{m}$ in \eqref{eq:monodomain} incorporates the conductivity properties of the tissue, and is defined as
\begin{equation}
    D_\mathrm{m} =
    \sigma_\mathrm{m}^\mathrm{l}\frac{F\mathbf f_0 \otimes F\mathbf f_0}{\|F\mathbf f_0\|^2} +
    \sigma_\mathrm{m}^\mathrm{t}\frac{F\mathbf s_0 \otimes F\mathbf s_0}{\|F\mathbf s_0\|^2} +
    \sigma_\mathrm{m}^\mathrm{n}\frac{F\mathbf n_0 \otimes F\mathbf n_0}{\|F\mathbf n_0\|^2}\;.
    \label{eq:conductivity}
\end{equation}
Here, $\sigma_\mathrm{m}^\mathrm{l}$, $\sigma_\mathrm{m}^\mathrm{t}$ and $\sigma_\mathrm{m}^\mathrm{n}$ are conductivities in the direction of sheets, fibers and normal to sheets respectively.

The formulations of both the monodomain equation \eqref{eq:monodomain} and the conductivity tensor \eqref{eq:conductivity} incorporate the so-called geometry-mediated mechano-electric feedback mechanisms, through the terms $J$ and $F$ \cite{collet2015numerical,salvador2021electromechanical,salvador2022role,trayanova2011whole} that account for the fact that the electrical stimulus is propagating in a deforming medium.

\begin{figure}
    \centering
    \begin{subfigure}{0.45\textwidth}
        \centering
        \raisebox{-0.5\height}{\includegraphics[width=0.45\textwidth]{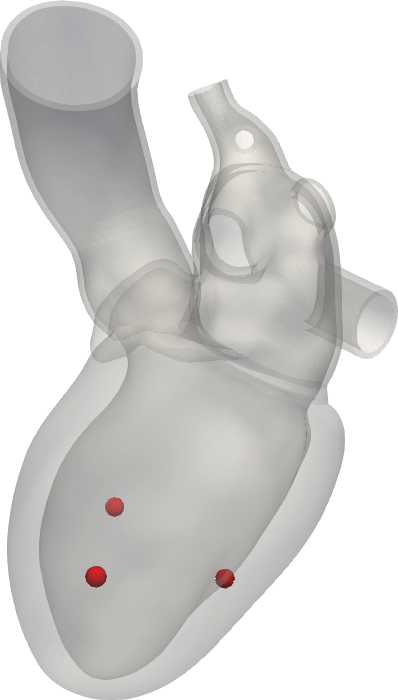}}
        \raisebox{-0.5\height}{\includegraphics[width=0.45\textwidth]{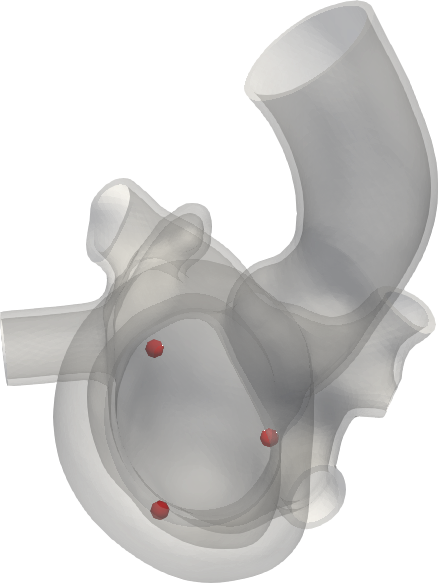}}
        \caption{}
        \label{fig:stimulus}
    \end{subfigure}
    \begin{subfigure}{0.45\textwidth}
        \centering
        \includegraphics[width=\textwidth]{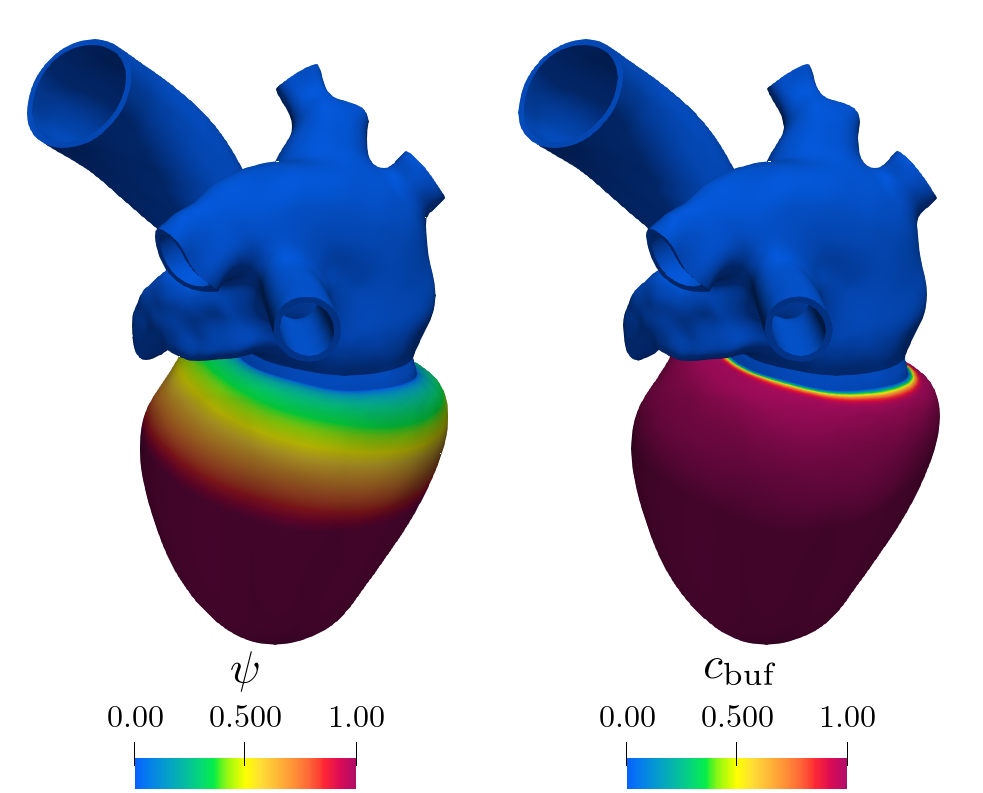}
        \caption{}
        \label{fig:buffer}
    \end{subfigure}
    \caption{(a) Stimulation points (red) on the left ventricular endocardium. (b) Functions $\psi$ (left) and $c_\mathrm{buf}$ used in the regularization of the interface between ventricle and atrioventricular ring.}
\end{figure}

Finally, $I_\mathrm{app}$ is a time-dependent forcing term that provides the initial stimulus. We impose an applied current on three points $\hat{\mathbf x}_\text{app}^0$, $\hat{\mathbf x}_\text{app}^1$ and $\hat{\mathbf x}_\text{app}^2$ on the endocardial surface of the ventricle, to trigger the electrical activation (\cref{fig:stimulus}). The applied current has the following analytical expression:
\begin{equation*}
    I_\mathrm{app}(\hat{\mathbf{x}}, t) = \begin{cases}
        \sum_{i = 0}^3 A_\mathrm{app}\,\exp{-\left(\frac{\|\hat{\mathbf x} - \hat{\mathbf x}_\mathrm{app}^i\|}{\sigma_\mathrm{app}}\right)^2} & \text{if } t \in (0, T_\mathrm{app}]\;, \\
        0 & \text{if } t > T_\mathrm{app}\;.
    \end{cases}
\end{equation*}
This mimics the effect of the Purkinje network \cite{costabal2016generating,del2022fast,landajuela2018numerical,romero2010effects,vergara2016coupled}, which is not included in our model. The stimulus is repeated every $T_\mathrm{hb} = \SI{800}{\milli\second}$ to obtain multiple heartbeats.

\subsection{Force generation model}

In response to the electrical excitation, muscular cells shorten, generating an active contractile force. As done for the electrophysiology, we model the ventricle $\hat\Omega\solid\lv$ as actively contracting, whereas all other subdomains are treated as mechanically passive.

We use the activation model RDQ20-MF presented in \cite{regazzoni2020biophysically}. The model is biophysically detailed, in the sense that it provides an explicit representation of the subcellular mechanisms leading to the generation of contractile force. Moreover, it includes the feedback between force generation and sarcomere length, responsible for the Frank-Starling mechanism \cite{jacob1992functional,katz2010physiology,klabunde2011cardiovascular,opie2004heart}, and between force generation and fiber shortening velocity \cite{bers2001excitation,katz2010physiology}. Both were found to be fundamental in capturing accurately the heart function, especially for what concerns the hemodynamics \cite{fedele2022comprehensive}.

The RDQ20-MF model is expressed in terms of a system of ODEs defined at each point of $\hat\Omega\solid\lv$:
\begin{equation}
    \begin{dcases}
        \pdv{\mathbf s}{t} = \mathbf F_\mathrm{act}\left(\mathbf s, [\mathrm{Ca}^{2+}]_i, SL, \pdv{SL}{t}\right) & \text{in } \hat\Omega\solid\lv \times (0, T)\;, \\
        \mathbf s = \mathbf s_0 & \text{in } \hat\Omega\solid\lv \times \{0\}\;.
    \end{dcases}
    \label{eq:activation}
\end{equation}
In the above system, $\mathbf s$ is a vector of variables defining the contraction state of cardiac cells, and $SL$ is the sarcomere length, defined as $SL = SL_0\sqrt{I_{4\mathrm f}}$. $SL_0$ is the sarcomere length at rest \cite{regazzoni2020biophysically} and $I_{4\mathrm f} = F\mathbf f_0 \cdot F\mathbf f_0$ is a measure of the stretch along the fibers.

In practice, to compute $SL$, we solve the following regularization problem:
\begin{equation}
    \begin{dcases}
        -\delta^2_{SL}\laplacian SL + SL = SL_0\sqrt{I_{4\mathrm f}} & \text{in } \hat\Omega\solid \times (0, T)\;, \\
        \delta^2_{SL}\grad SL\cdot\mathbf n = 0 & \text{on } \partial\hat\Omega\solid \times (0, T)\;,
    \end{dcases}
    \label{eq:sl}
\end{equation}
with $\delta_{SL}$ a suitable regularization radius parameter. This has the effect of preventing sharp variations of $SL$ over spatial scales smaller than $\delta_{SL}$ \cite{regazzoni2022cardiac}.

The generated active force is computed at every point $\hat{\mathbf x} \in \hat\Omega\solid\lv$ as a function $T_\mathrm{act}(\mathbf s)$ of the contraction state $\mathbf s$. For the precise definition of $\mathbf F_\mathrm{act}$ and $T_\mathrm{act}$, we refer to \cite{regazzoni2020biophysically}. The force-velocity relationship may yield instabilities, which are prevented with a numerically consistent stabilization as described in \cite{regazzoni2021oscillation}.

The initial state $\mathbf s_0$ is obtained by solving the system \eqref{eq:activation} uncoupled from the other models, until a steady state solution is reached.

\subsection{Solid mechanics}
\label{sec:mechanics}

The deformation of the heart wall is modeled in a Lagrangian reference framework with the elastodynamics equation \cite{ogden2013non}. The vector $\mathbf d$ describes the displacement of the solid domain, so that
\begin{equation*}
    \mathcal L\solid(\hat{\mathbf x}, t) = \hat{\mathbf x} + \mathbf d(\hat{\mathbf x}, t)\qquad \hat{\mathbf x} \in \hat\Omega\solid\;, t \in (0, T)\;.
\end{equation*}
Let $\hat\Gamma\solid^\mathrm{epi} \subset \partial\hat\Omega\solid\backslash\hat\Sigma$ be the epicardial boundary, corresponding to the outer wall of the heart, $\hat\Gamma\solid^\mathrm{PV} \subset \partial\hat\Omega\solid\backslash\hat\Sigma$ be the portion of the boundary corresponding to the pulmonary veins sections, and $\hat\Gamma\solid^\mathrm{ao} \subset \partial\hat\Omega\solid\backslash\hat\Sigma$ be the one corresponding to the AA terminal section (\cref{fig:domain-boundary}). The equations for solid mechanics read
\begin{equation}
    \begin{dcases}
        \rho\solid\pdv[2]{\mathbf d}{t} - \div P\solid(\mathbf d, \mathbf s) = \mathbf 0 & \text{in } \hat\Omega\solid \times (0, T)\;, \\
        \mathbf d = \mathbf 0 & \text{on } \left(\hat\Gamma\solid^\mathrm{PV} \cup \hat\Gamma\solid^\mathrm{ao}\right) \times (0, T)\;, \\
        \begin{multlined} P\solid(\mathbf d, \mathbf s)\mathbf n = -(\mathbf n \otimes \mathbf n)\left(K_\perp^\mathrm{epi}\mathbf d + C_\perp^\mathrm{epi}\pdv{\mathbf d}{t}\right) \\ - (I - \mathbf n \otimes \mathbf n)\left(K_\parallel^\mathrm{epi}\mathbf d + C_\parallel^\mathrm{epi}\pdv{\mathbf d}{t}\right)\end{multlined} & \text{on } \hat\Gamma\solid^\mathrm{epi} \times (0, T)\;, \\
        \mathbf d = \mathbf d_0 & \text{in } \hat\Omega\solid \times \{0\}\;.
    \end{dcases}
    \label{eq:elastodynamics}
\end{equation}

In the above system, $\rho\solid$ is the density of the solid and $P\solid(\mathbf d, \mathbf s)$ is the first Piola-Kirchhoff stress tensor. The stress tensor incorporates both the passive and the active mechanical properties of the material. We work in the active stress framework \cite{ambrosi2012active,goktepe2010electromechanics,gurev2011models,pathmanathan2010cardiac,quarteroni2019mathematical,regazzoni2022cardiac,smith2004multiscale}, decomposing the stress tensor into the sum of a passive and an active term:
\begin{equation*}
    P\solid(\mathbf d, \mathbf s) = P_\mathrm{pas}(\mathbf d) + P_\mathrm{act}(\mathbf d, \mathbf s)\;.
\end{equation*}

The passive part is obtained as the derivative of a strain energy function $\mathcal{W}$ that characterizes the mechanical properties of the material \cite{ogden2013non}:
\begin{equation*}
    P_\mathrm{pas}(\mathbf d) = \pdv{\mathcal{W}}{F}\;.
\end{equation*}
We use the Guccione anisotropic constitutive law \cite{guccione1991finite,usyk2002computational} in $\hat\Omega\solid\lv$ and in $\hat\Omega\solid\la$. In $\hat\Omega\solid\ao$ and $\hat\Omega\solid\ring$ we use an isotropic neo-Hookean material model \cite{ogden2013non}. In both cases, we treat the solid as nearly incompressible, by penalizing local volume variations in the strain energy function. See \cref{app:constitutive-laws} for the definitions of the constitutive laws.

The active part of the stress tensor is computed as a function of the activation state, through
\begin{equation*}
    P_\mathrm{act}(\mathbf d, \mathbf s) = T_\mathrm{act}(\mathbf s)\frac{F\mathbf f_0 \otimes \mathbf f_0}{\sqrt{I_{4\mathrm f}}}\;.
\end{equation*}

The boundary conditions imposed on the mechanics equations are a homogeneous Dirichlet condition on $\hat\Gamma\solid^\mathrm{PV}$ and $\hat\Gamma\solid^\mathrm{ao}$ \cite{gerach2021electro}, and a generalized visco-elastic Robin boundary condition on $\hat\Gamma\solid^\mathrm{epi}$. The latter has the purpose of modeling the presence of the pericardial sac, that surrounds the heart and provides mechanical support \cite{pfaller2019importance,strocchi2020simulating}. Finally, no boundary conditions are prescribed on $\hat\Sigma$ since the interface conditions for the fluid-solid coupling are set on that portion of $\partial\hat\Omega\solid$ (see \cref{sec:fsi}).

\subsubsection{Initial conditions}

We start from a model of the left heart that corresponds to a configuration reconstructed from CT and MRI scans \cite{zygote}. Such configuration is subject to the load of blood pressure on the endocardial wall, while the elastodynamics equations \eqref{eq:elastodynamics} are formulated assuming that $\hat\Omega\solid$ is stress-free. Therefore, as a preprocessing step, we compute a stress-free reference configuration using the unloading algorithm presented in \cite{regazzoni2022cardiac}.

From there, we determine the initial condition for the solid displacement $\mathbf d_0$ by solving a quasi-static ramp, as described in \cite{regazzoni2022cardiac}, imposing an initial pressure $p_0$ on the endocardial wall. We choose $p_0$ to match typical values at the end of the diastolic phase. We remark that different initial pressures $p_0\la$, $p_0\ring$, $p_0\lv$ and $p_0\ao$ are prescribed on the boundaries of the different chambers \cite{fedele2022comprehensive}.

\subsubsection{Material interface regularization}

The different material models used result in sharp discontinuities across the interface between the atrioventricular ring subdomain and the atrium and ventricle subdomains. These do not constitute an issue from the point of view of solid mechanics. However, they lead to the formation of corners at the fluid-solid interface, which pose a relevant numerical issue for the fluid domain displacement inducing mesh elements distortion (see \cref{sec:fluid_domain}). To overcome it, we regularize the interface between $\hat\Omega\solid\lv$ and $\hat\Omega\solid\ring$ as follows. As a preprocessing step, a Laplace problem is solved to obtain a function $\psi$ that varies smoothly from \num{0} in $\hat\Omega\solid\ring$ to \num{1} in the apical portion of $\hat\Omega\solid\lv$. We then set
\begin{equation*}
    c_\mathrm{buf}(\hat{\mathbf x}) = \frac{1}{2}\left[1 - \cos\left(\frac{\pi \min\{\psi, \psi_\mathrm{th}\}}{\psi_\mathrm{th}}\right)\right]\;,
\end{equation*}
where $\psi_\mathrm{th}$ is a threshold value that controls the size of the regularization region. Both $\psi$ and $c_\mathrm{buf}$ are depicted in \cref{fig:buffer}. Finally, we redefine the stress tensor in $\hat\Omega\solid\lv$ as a convex combination of the Guccione (G) and neo-Hookean (NH) stress tensors so that the transition between the material of the ventricle (or the atrium) and the atrioventricular ring is smooth:
\begin{equation*}
    P_\mathrm{pas}(\mathbf d) = c_\mathrm{buf}\,P_\mathrm{G}(\mathbf d) + (1 - c_\mathrm{buf})\,P_\mathrm{NH}(\mathbf d)\;.
\end{equation*}
The regularization is included at the interfaces between $\hat\Omega\solid\lv$ and $\hat\Omega\solid\ring$ and between $\hat\Omega\solid\la$ and $\hat\Omega\solid\ring$.

\subsection{Fluid domain displacement}
\label{sec:fluid_domain}

We account for the motion of $\Omega\fluid$ in the arbitrary Lagrangian-Eulerian (ALE) framework \cite{donea1982arbitrary,hughes1981lagrangian,stein2003mesh}. The vector $\mathbf d\ale$ represents the displacement of the fluid domain, such that
\begin{equation*}
    \mathcal{L}\fluid(\hat{\mathbf x}, t) = \hat{\mathbf x} + \mathbf d\ale(\hat{\mathbf x}, t)\qquad \hat{\mathbf x} \in \hat\Omega\fluid\;, t \in (0, T)\;.
\end{equation*}
Inspired by \cite{hoffman2011unified,spuhler2021interface}, we use a fictitious non-linear solid material to model the displacement of the fluid domain, so that $\mathbf d\ale$ solves at every time $t$ the following stationary problem:
\begin{equation}
    \label{eq:lifting}
    \begin{dcases}
        -\div P\ale(\mathbf d\ale) = \mathbf 0 & \text{in }\hat\Omega\fluid\;, \\
        \mathbf d\ale = \mathbf d & \text{on }\hat\Sigma\;, \\
        \mathbf d\ale = \mathbf 0 & \text{on }\hat\Gamma\fluid^\mathrm{PV} \cup \hat\Gamma\fluid^\mathrm{ao}\;,
    \end{dcases}
\end{equation}
wherein
\begin{gather*}
    P\ale = \frac{1}{q}\left(I - \left(F\ale F\ale^T\right)^{-1}\right)\;, \\
    F\ale = I + \grad\mathbf d\ale\;.
\end{gather*}
In the above, $q$ is a scale-invariant mesh quality metric that has the purpose of stiffening the regions of the fluid domain with highly distorted element, aiming at preventing solver breakdown. It is defined element-wise as
\begin{equation*}
    q(\hat{\mathbf x}) = \frac{|D\ale|^2_\mathrm{F}}{3(\det D\ale)^\frac{2}{3}}\;,
\end{equation*}
where $D\ale = F\ale\,\grad\mathcal{M}$ and $\mathcal{M}$ is the linear mapping from the unit simplex to the element in current configuration.

We remark that \eqref{eq:lifting} is non-linear, since $P\ale$ is a non-linear function of $\mathbf d\ale$. Therefore, its solution is significantly more complex and costly than that of other typically used operators, such as the harmonic extension operator \cite{bucelli2022partitioned,zhang2001analysis,zingaro2022geometric} or linear elasticity \cite{johnson1994mesh,stein2003mesh}. On the other hand, however, we found that this results in an increased robustness with respect to large deformations, preventing the inversion of mesh elements, that would result in the breakdown of the numerical solver.

We define the ALE velocity as the time derivative of $\mathbf d\ale$, pushed forward to the current configuration, namely
\begin{equation*}
    \mathbf u\ale(\mathbf x, t) = {\pdv{\mathbf d\ale}{t}}\left(\mathcal{L}\fluid^{-1}(\mathbf x, t), t\right)\;.
\end{equation*}

\subsection{Fluid dynamics}

The evolution of fluid velocity $\mathbf u$ and pressure $p$ is prescribed by the Navier-Stokes equations for a Newtonian, incompressible fluid \cite{quarteroni2017numerical}:
\begin{subnumcases}{\label{eq:ns}}
    \rho\fluid\left[\pdv{\mathbf u}{t} + \left(\left(\mathbf u - \mathbf u\ale\right) \cdot\grad\right)\mathbf u\right] - \div\sigma\fluid(\mathbf u, p) + \boldsymbol{\mathcal{R}}(\mathbf u, \mathbf u\ale) = \mathbf 0 & in $\Omega\fluid \times (0, T)$, \label{eq:nsmomentum} \\
    \div \mathbf u = 0 &  in $\Omega\fluid \times (0, T)$, \\
    \sigma\fluid(\mathbf u, p)\mathbf n = -p_\mathrm{in}(t)\mathbf n & on $ \Gamma\fluid^\mathrm{PV} \times (0, T)$, \label{eq:nsinlet} \\
    \sigma\fluid(\mathbf u, p)\mathbf n = -p_\mathrm{out}(t)\mathbf n & on $ \Gamma\fluid^\mathrm{ao} \times (0, T)$, \label{eq:nsoutlet} \\
    \mathbf u = \mathbf 0 & in $ \Omega\fluid \times \{0\}$.
\end{subnumcases}
In the above system, $\rho\fluid$ denotes the fluid density, and the stress tensor $\sigma\fluid$ is given by
\begin{equation*}
    \sigma\fluid(\mathbf u, p) = \mu\fluid \left(\grad\mathbf u + \grad\mathbf u^T\right) - p I\;,
\end{equation*}
where $\mu\fluid$ is the dynamic viscosity. $\boldsymbol{\mathcal{R}}(\mathbf u)$ is a resistive term that accounts for the presence of valves in a penalty-based approach (see \cref{sec:valves}) \cite{corti2022impact,fedele2017patient,fumagalli2020image,fumagalli2022image,fumagalli2021reduced,zingaro2022geometric}.

In analogy with the solid mechanics problem (\cref{sec:mechanics}), no boundary conditions are prescribed on $\Sigma$ since, on that portion of $\partial\Omega\fluid$, the interface conditions for the fluid-solid coupling are imposed (see \cref{sec:fsi}). The functions $p_\mathrm{in}(t)$ and $p_\mathrm{out}(t)$ are pressures provided by the circulation model (see \cref{sec:circulation}) for the pulmonary venous and systemic arterial circulation compartments, respectively.

\subsubsection{Valve modeling}
\label{sec:valves}

Heart compartments are separated by valves that prevent reverse flow \cite{katz2010physiology}. For a left heart model, the mitral valve (MV) separates the LV and the LA, and the aortic valve (AV) separates the LV and the AA. Accurate modeling and numerical simulation of the valves is challenging: the valves are thin structures, they undergo large displacements, and contact phenomena play a major role in their physiological function \cite{astorino2009fluid,dabiri2019tricuspid,einstein2010fluid,hirschhorn2020fluid,hsu2014fluid,luraghi2017evaluation,spuhler20183d,terahara2020heart}. Our interest mainly lies in the macroscopic effect of the opening and closing of the valves onto the blood flow, rather than on an accurate description of valves motion. In particular, we want our model to describe the role of valves in ensuring the correct direction of the flow through the heart, and to capture the formation of jets and vortices associated to the presence of valve leaflets.

Therefore, we use a reduced approach for the modeling of valves, based on the Resistive Immersed Implicit Surface (RIIS) method \cite{fedele2017patient,fumagalli2020image,fumagalli2022image,zingaro2022geometric,zingaro2022modeling}. Each valve, at any time $t \in (0, T)$, is represented by a surface $\Gamma_\mathrm{k}^t$, $\mathrm{k} \in \{\mathrm{MV}, \mathrm{AV}\}$, immersed in the fluid domain. For $\mathrm{k} \in \{\mathrm{MV}, \mathrm{AV}\}$, let $\phi_\mathrm{k}^t(\mathbf x)$ be the signed distance function from the surface $\Gamma_\mathrm{k}^t$. The Navier-Stokes momentum equation \eqref{eq:nsmomentum} includes the penalization term $\boldsymbol{\mathcal{R}}(\mathbf u, \mathbf u\ale)$ forcing the fluid velocity to match the valve velocity (with respect to the moving frame of the domain) near the immersed surfaces:
\begin{equation*}
    \boldsymbol{\mathcal{R}}(\mathbf u, \mathbf u\ale) = \sum_{\mathrm{k} \in \{\mathrm{MV}, \mathrm{AV}\}} \frac{R_\mathrm{k}}{\epsilon_\mathrm{k}}\delta_{\epsilon_\mathrm{k}}(\phi_\mathrm{k}^t(\mathbf x))(\mathbf u - \mathbf u\ale - \mathbf u_{\Gamma_\mathrm{k}})\;.
\end{equation*}
In the above, $R_\mathrm{k}$ are resistance parameters, $\epsilon_\mathrm{k}$ are valves half-thicknesses, $\mathbf u_{\Gamma_\mathrm{k}}$ is the valve velocity with respect to the moving domain, and $\delta_{\epsilon_\mathrm{k}}$ is a smoothed delta function, defined as
\begin{equation*}
    \delta_{\epsilon_\mathrm{k}}(y) = \begin{dcases}
        \frac{1}{2\epsilon_\mathrm{k}}\left(1 + \cos\left(\frac{\pi y}{\epsilon_\mathrm{k}}\right)\right) & \text{if } |y| \leq \epsilon_\mathrm{k}\;, \\
        0 & \text{if } |y| > \epsilon_\mathrm{k}\;.
    \end{dcases}
\end{equation*}

We account for the opening and closing of valves by deforming the corresponding surfaces. Although we do not include any FSI model for the valves, their opening and closing times are determined based on the pressure jump across the immersed surfaces while the durations of opening and closing are prescribed. We refer to \cite{zingaro2022geometric} and \cref{app:valves} for more details.

\subsubsection{Flow stabilization and turbulence modeling}

Blood flow in the heart chambers is characterized by a regime of transition to turbulence \cite{bluestein1994transition,verkaik2012coupled,vignon2010outflow,zingaro2021hemodynamics,zingaro2022geometric}. To account for that, we use the VMS-LES model for fine scales \cite{bazilevs2007variational,zingaro2021hemodynamics}. We refer the interested reader to \cite{bazilevs2007variational} for details on the VMS-LES approach to turbulence modeling.

We set Neumann boundary conditions on boundaries $\Gamma\fluid^\mathrm{PV}$ and $\Gamma\fluid^\mathrm{ao}$. Neumann conditions may cause instability phenomena in case of inflow \cite{bertoglio2014tangential,moghadam2011comparison}. Therefore, we make use of backflow stabilization in the form of inertial stabilization as presented in \cite{bertoglio2014tangential}: Neumann boundary conditions are modified by imposing
\begin{equation*}
    \begin{dcases}
        \sigma\fluid(\mathbf u, p)\mathbf n = -p_\mathrm{in}(t)\mathbf n + \beta\frac{\rho\fluid}{2}|\mathbf u\cdot\mathbf n|_-\mathbf u & \text{on } \Gamma\fluid^\mathrm{PV} \times (0, T)\;, \\
        \sigma\fluid(\mathbf u, p)\mathbf n = -p_\mathrm{out}(t)\mathbf n + \beta\frac{\rho\fluid}{2}|\mathbf u\cdot\mathbf n|_-\mathbf u & \text{on } \Gamma\fluid^\mathrm{ao} \times (0, T)\;, \\
    \end{dcases}
\end{equation*}
where $|\mathbf u\cdot\mathbf n|_- = \min\{\mathbf u\cdot\mathbf n, 0\}$ and $\beta = \num{1}$.

\subsection{Fluid-structure interaction}
\label{sec:fsi}

The fluid and solid models are coupled by kinematic and dynamic interface conditions on $\Sigma$ that prescribe the continuity of velocity and of stresses \cite{bazilevs2013computational,bucelli2022partitioned}:
\begin{equation*}
    \begin{dcases}
        \mathbf u = \pdv{\mathbf d}{t} & \text{on }\Sigma\times (0, T)\;, \\
        \sigma\fluid(\mathbf u, p)\mathbf n = \sigma\solid(\mathbf d, \mathbf s)\mathbf n & \text{on }\Sigma\times (0, T)\;.
    \end{dcases}
\end{equation*}
In the above, $\sigma\solid(\mathbf d)$ is the Cauchy stress tensor for the solid, related to the first Piola-Kirchhoff tensor by
\begin{equation*}
    J \sigma\solid(\mathbf d, \mathbf s) = F P\solid(\mathbf d, \mathbf s)^T\;.
\end{equation*}

\begin{figure}
    \centering
    \includegraphics[width=0.9\textwidth]{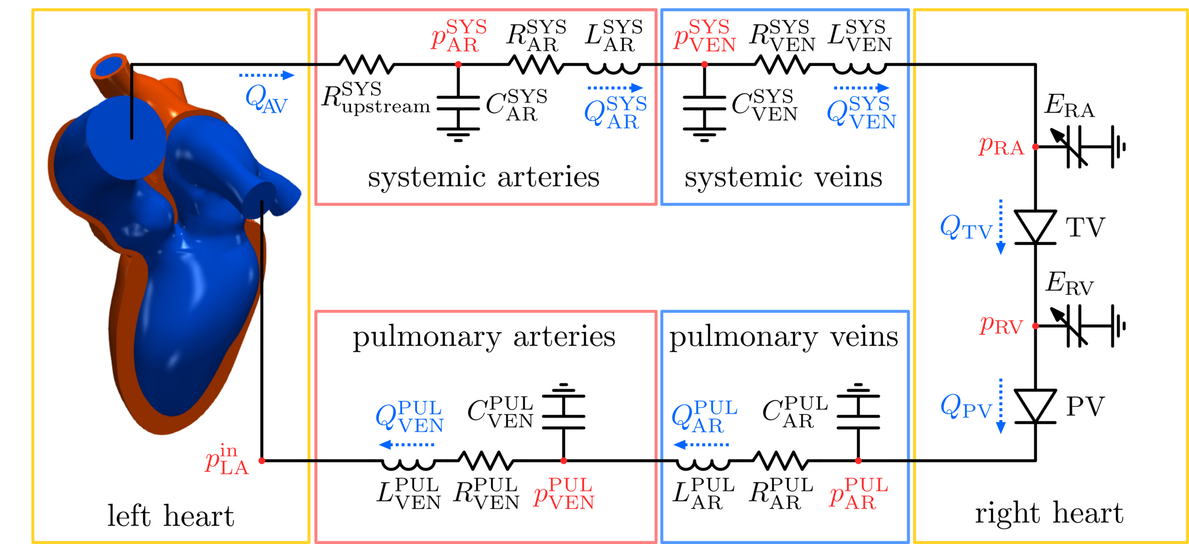}
    \caption{Circuital (0D) representation for the circulation model coupled to the EMF model. Refer to \cref{app:circulation} for the model equations.}
    \label{fig:circulation}
\end{figure}

\subsection{Circulation model}
\label{sec:circulation}

Inlet and outlet conditions for the fluid model are provided by a zero-dimensional circulation model \cite{blanco20103d,hirschvogel2017monolithic,regazzoni2022cardiac,zingaro2022geometric}, representing the whole closed-loop circulation with a lumped-parameter approach. It includes four different circulation compartments (systemic arterial, systemic venous, pulmonary arterial and pulmonary venous) as well as the four heart chambers. Using an electric circuit analogy \cite{hirschvogel2017monolithic} represented in \cref{fig:circulation}, the circulation compartments are modeled as RLC circuits, heart chambers are modeled as time-varying elastances, and valves are modeled as non-ideal diodes.

In order to couple it to the three-dimensional FSI model, we remove from the full circulation model those compartments that have a three-dimensional description, that is the left atrium, the left ventricle, and mitral and aortic valves \cite{zingaro2022geometric}. The remaining unknowns are the volumes and pressures of the right atrium and ventricle ($V_\mathrm{RA}$, $p_\mathrm{RA}$, $V_\mathrm{RV}$ and $p_\mathrm{RV}$), pressures and flow rates through the four circulation compartments, except for the pulmonary venous flow rate ($p\ar\sys$, $Q\ar\sys$, $p\ven\sys$, $Q\ven\sys$, $p\ar\pul$, $Q\ar\pul$ and $p\ven\pul$), and flow rates through tricuspid and pulmonary valves ($Q_\mathrm{PV}$ and $Q_\mathrm{TV}$). Collecting all the circulation variables into the vector $\mathbf c(t)$, the circulation problem can be stated in compact form as a differential-algebraic equation:
\begin{equation*}
    \mathbf F_\mathrm{circ}\left(\pdv{\mathbf c}{t}, \mathbf c, t\right) = \mathbf 0\;.
    \label{eq:circulation}
\end{equation*}
Refer to \cref{app:circulation} for the definition of $\mathbf F_\mathrm{circ}$.

Following \cite{quarteroni2016geometric,zingaro2022geometric}, the 0D circulation model is coupled to the 3D fluid model by imposing the continuity of stresses and of fluxes at the 3D-0D interface: for $t \in (0, T)$,
\begin{align*}
    p_\mathrm{in}(t) &= p^\text{in}_\text{LA}(t)\;, \\
    p_\mathrm{out}(t) &= p\ar\sys(t) + R\sys_\mathrm{upstream}Q_\mathrm{AV}(t)\;, \\
    Q\ven\pul(t) &= -\int_{\Gamma\fluid^\mathrm{PV}}(\mathbf u - \mathbf u\ale)\cdot\mathbf n d\gamma\;, \\
    Q_\mathrm{AV}(t) &= \int_{\Gamma\fluid^\mathrm{ao}}(\mathbf u - \mathbf u\ale)\cdot\mathbf n d\gamma\;.
\end{align*}
We remark that, with respect to \cite{zingaro2022geometric}, we include the resistance term $R\sys_\mathrm{upstream}Q_\mathrm{AV}(t)$ when coupling the outlet to the systemic arterial compartment. This has the role of avoiding reflections of pressure waves at the outlet, that would otherwise result in unphysical oscillations \cite{janela2010comparing}. $p^\text{in}_\text{LA}(t)$ is the pressure downstream of the pulmonary veins compartment (\cref{fig:circulation}).

\begin{figure}
    \centering
    \includegraphics[width=0.85\textwidth]{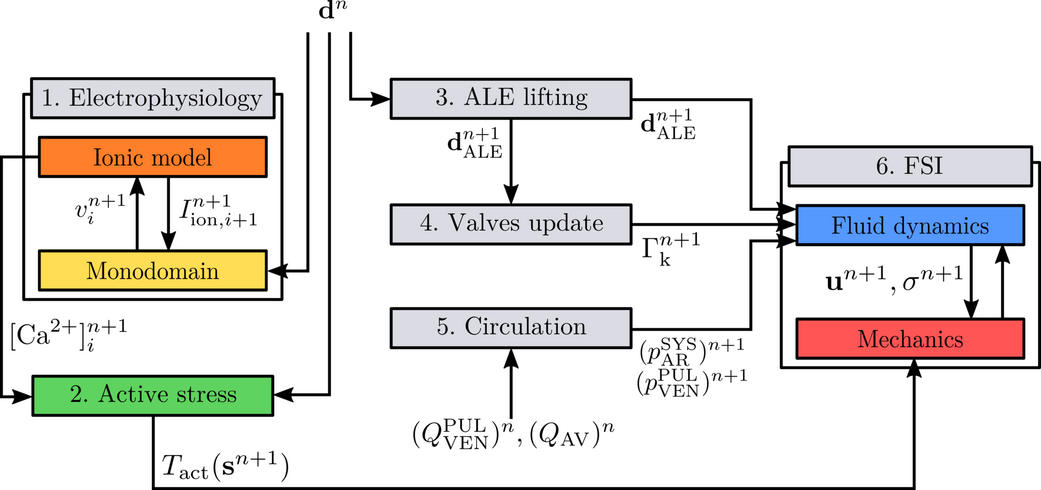}
    \caption{Time advancing scheme for the coupled model. The numbers correspond to the steps described in \cref{sec:time-discretization}.}
    \label{fig:timescheme}
\end{figure}

\section{Numerical discretization}
\label{sec:numerical}

We detail the discretization strategy used to numerically solve the fully coupled EMF problem, starting from the semi-discretization in time (\cref{sec:time-discretization}) and then describing spatial discretization (\cref{sec:space-discretization}).

\subsection{Time discretization}
\label{sec:time-discretization}

We introduce a partition of the time interval $(0, T)$ into sub-intervals $(t_i, t_{i+1}]$, with $i = 0, 1, \dots, N_T$, $t_0 = 0$ and $t_{N_T} = T$, and such that $t_{i+1} - t_i = \Delta t$ for all $i$. From here on, we denote the approximation of any of the solution variables at a given timestep $t_n$ with the superscript $n$, e.g. $\mathbf u^n \approx \mathbf u(t_n)$.

We employ a staggered scheme for the coupling of the different systems of equations, schematically represented in \cref{fig:timescheme}. The scheme is based on solving separately the different subproblems, using an explicit coupling whenever stability is not a concern. We remark that we use an implicit coupling of the FSI problem \cite{bucelli2022partitioned}. The scheme is derived from that in \cite{regazzoni2022cardiac} for the electromechanics of the left ventricle, adapted to include the three dimensional description of the blood hemodynamics.

We use finite differences for the approximation of time derivatives appearing in the different model equations. Electrophysiology is characterized by a faster dynamics than the other models, and requires a smaller timestep for an accurate solution \cite{regazzoni2022cardiac}. Therefore, we discretize it with a timestep $\Delta t_\mathrm{EP} = \frac{\Delta t}{N_\mathrm{EP}}$, with $N_\mathrm{EP} \in \mathbb{N}$. Given the solution variables and the domain $\Omega\fluid^n$ at timestep $t_n$, in order to compute the solution at timestep $t_{n+1}$ we perform the following steps:
\begin{enumerate}
    \item Electrophysiology: we solve $N_\mathrm{EP}$ time advancing steps of the electrophysiology problem, that is: setting $v^{n+1}_0 = v^n$ and $\mathbf w^{n+1}_0 = \mathbf w^n$, for $i = 0, 1, \dots, N_\mathrm{EP}-1$:
    \begin{enumerate}
        \item we compute $\mathbf w^{n+1}_{i+1}$ and $\mathbf z^{n+1}_{i+1}$ by solving the ionic model
        \begin{equation}
            \left\{\begin{aligned}
                \frac{\mathbf w^{n+1}_{i+1} - \mathbf w^{n+1}_i}{\Delta t_\mathrm{EP}} &= \mathbf F_\mathrm{ion}^\mathbf{w}(v^{n+1}_i, \mathbf w^{n+1}_{i+1}) & \text{in } \hat\Omega\solid\lv\;, \\
                \frac{\mathbf z^{n+1}_{i+1} - \mathbf z^{n+1}_i}{\Delta t_\mathrm{EP}} &= \mathbf F_\mathrm{ion}^\mathbf{z}(v^{n+1}_i, \mathbf w^{n+1}_i, \mathbf z^{n+1}_i) & \qquad\text{in } \hat\Omega\solid\lv\;.
            \end{aligned}\right.
            \label{eq:ionic-time-discrete}
        \end{equation}
        We use an implicit-explicit (IMEX) scheme, with an implicit discretization of gating variables and an explicit discretization of ionic concentrations \cite{regazzoni2020mathematical,regazzoni2022cardiac}.
        \item Then, we solve the monodomain equation to compute $v^{n+1}_{i+1}$:
        \begin{equation}
            \begin{dcases}
                \begin{multlined}
                J^n \chi C_\text{m} \frac{v^{n+1}_{i+1} - v^{n+1}_i}{\Delta t_\mathrm{EP}} - \div(J^n (F^n)^{-1} D_\mathrm{m} (F^n)^{-T} \grad v^{n+1}_{i+1}) \\ + J^n \chi I_\mathrm{ion}(v^{n+1}_i, \mathbf w^{n+1}_{i+1}, \mathbf z^{n+1}_{i+1}) = J^n \chi I_\mathrm{app}(\hat{\mathbf x}, t^n + (i + 1)\Delta t_\mathrm{EP})\end{multlined} & \text{in } \hat\Omega\solid\lv\;, \\[1em]
                J^n (F^n)^{-1}D_\mathrm{m}(F^n)^{-T}\grad v^{n+1}_{i+1}\cdot\mathbf n = 0 & \text{on } \partial\hat\Omega\solid\lv\;.
            \end{dcases}
            \label{eq:monodomain-time-discrete}
        \end{equation}
        System \eqref{eq:monodomain-time-discrete} relies on a semi-implicit discretization, since the ionic current term is computed using the tramsmembrane potential at previous subiteration, $v^{n+1}_i$. Therefore, the problem is linear.
    \end{enumerate}
    We set $v^{n+1} = v^{n+1}_{N_\mathrm{EP}}$, $\mathbf w^{n+1} = \mathbf w^{n+1}_{N_\mathrm{EP}}$ and $\mathbf z^{n+1} = \mathbf z^{n+1}_{N_\mathrm{EP}}$;

    \item Force generation model: we compute $SL^{n}$ solving \eqref{eq:sl} using $\mathbf d^n$ to compute $I_{4\text{f}}$, then we solve the activation model to compute $\mathbf s^{n+1}$:
    \begin{equation}
        \frac{\mathbf s^{n+1} - \mathbf s^n}{\Delta t} = \mathbf F_\mathrm{act}\left(s^n, [\text{Ca}^{2+}]_i^{n+1}, SL^n, \frac{SL^n - SL^{n-1}}{\Delta t}\right)\;;
        \label{eq:activation-time-discrete}
    \end{equation}

    \item Fluid domain displacement: we compute $\mathbf d\ale^{n+1}$ by solving
    \begin{equation}
        \begin{dcases}
            -\div P\ale(\mathbf d\ale^{n+1}) = \mathbf 0 & \text{in } \hat\Omega\fluid\;, \\
            \mathbf d\ale^{n+1} = \mathbf d^n  & \text{on } \hat\Sigma\;, \\
            \mathbf d\ale^{n+1} = \mathbf 0 & \text{on } \hat\Gamma\fluid^\mathrm{PV} \cup \hat\Gamma\fluid^\mathrm{ao}\;.
        \end{dcases}
        \label{eq:lifting-time-discrete}
    \end{equation}
    We remark that we compute $\mathbf d\ale^{n+1}$ from the displacement $\mathbf d^n$ from the previous timestep, so that the geometric FSI coupling condition is treated explicitly. We then update the fluid domain $\Omega\fluid^{n+1}$ according to the displacement $\mathbf d\ale^{n+1}$ and compute \[\mathbf u\ale^{n+1}(\mathbf x) = \frac{\mathbf d\ale^{n+1}\left(\left(\mathcal{L}\fluid^{n+1}\right)^{-1}(\mathbf x)\right) - \mathbf d\ale^n\left(\left(\mathcal{L}\fluid^{n}\right)^{-1}(\mathbf x)\right)}{\Delta t}\;;\]

    \item Valves position: we update the position of the valves according to their opening state and to $\mathbf d\ale^{n+1}$, computing the surfaces $\Gamma_\text{MV}^{n+1}$ and $\Gamma_\text{AV}^{n+1}$ (see \cref{app:valves});

    \item Circulation: we compute the flow rates at the 3D-0D interface as
    \begin{gather*}
        (Q\pul\ven)^{n+1} = -\int_{\Gamma\fluid^\mathrm{PV}}(\mathbf u^{n} - \mathbf u\ale^{n+1})\cdot\mathbf n d\gamma\;, \\
        Q_\mathrm{AV}^{n+1} = \int_{\Gamma\fluid^\mathrm{ao}}(\mathbf u^{n} - \mathbf u\ale^{n+1})\cdot\mathbf n d\gamma\;,
    \end{gather*}
    then advance the circulation model \eqref{eq:circulation} with an explicit Runge-Kutta scheme \cite{quarteroni2017numerical}, to compute $\mathbf c^{n+1}$;

    \item Fluid-structure interaction: we solve the FSI problem to compute $\mathbf d^{n+1}, \mathbf u^{n+1}, p^{n+1}$:
    \begin{equation}
        \begin{dcases}
            \rho\fluid\frac{\mathbf d^{n+1} - 2\mathbf d^n + \mathbf d^{n-1}}{\Delta t^2} - \div P\solid(\mathbf d^{n+1}, \mathbf s^{n+1}) = \mathbf 0 & \text{in } \hat\Omega\solid\;, \\
            \mathbf d^{n+1} = \mathbf 0 & \text{on } \hat\Gamma\solid^\mathrm{PV} \cup \hat\Gamma\solid^\mathrm{ao}\;, \\
            \begin{multlined}P\solid(\mathbf d^{n+1}, \mathbf s^{n+1})\mathbf n = -(\mathbf n \otimes \mathbf n)\left(K_\perp^\mathrm{epi} \mathbf d^{n+1} + C_\perp^\mathrm{epi} \frac{\mathbf d^{n+1} - \mathbf d^{n}}{\Delta t}\right) \\ - (I - \mathbf n \otimes \mathbf n)\left(K_\parallel^\mathrm{epi} \mathbf d^{n+1} + C_\parallel^\mathrm{epi} \frac{\mathbf d^{n+1} - \mathbf d^{n}}{\Delta t}\right)\end{multlined} & \text{on } \hat\Gamma\solid^\mathrm{epi}\;, \\[1em]
            \begin{multlined}\rho\fluid\left[\frac{\mathbf u^{n+1} - \mathbf u^n}{\Delta t} + ((\mathbf u^{n+1} - \mathbf u\ale^{n+1})\cdot\grad)\mathbf u^{n+1}\right] \\ - \div\sigma\fluid(\mathbf u^{n+1}, p^{n+1}) + \boldsymbol{\mathcal{R}}(\mathbf u^{n+1}, \mathbf u\ale^{n+1}) = \mathbf 0\end{multlined} & \text{in } \Omega\fluid\;, \\[1em]
            \div\mathbf u^{n+1} = 0 & \text{in } \Omega\fluid\;, \\
            \sigma\fluid(\mathbf u^{n+1}, p^{n+1})\mathbf n = -\left(p\pul\ven\right)^{n+1}\mathbf n + \rho\fluid\frac{\beta}{2}|\mathbf u^{n+1}\cdot\mathbf n|_- \mathbf u^{n+1} & \text{on } \Gamma\fluid^\mathrm{PV}\;, \\
            \sigma\fluid(\mathbf u^{n+1}, p^{n+1})\mathbf n = -\left(p\ar\sys\right)^{n+1}\mathbf n + \rho\fluid\frac{\beta}{2}|\mathbf u^{n+1}\cdot\mathbf n|_- \mathbf u^{n+1} & \text{on } \Gamma\fluid^\mathrm{ao}\;, \\
            \mathbf u^{n+1} = \frac{\mathbf d^{n+1} - \mathbf d^n}{\Delta t} & \text{on } \Sigma\;, \\
            \sigma\fluid(\mathbf u^{n+1}, p^{n+1})\mathbf n = \sigma\solid(\mathbf d^{n+1}, \mathbf s^{n+1})\mathbf n & \text{on } \Sigma\;.
        \end{dcases}
        \label{eq:fsi-time-discrete}
    \end{equation}
\end{enumerate}

For our numerical simulations, we set $\Delta t = \SI{0.2}{\milli\second}$ and $N_\text{EP} = \num{2}$, so that $\Delta t_\text{EP} = \SI{0.1}{\milli\second}$.

\subsection{Space discretization}
\label{sec:space-discretization}

\begin{figure}
    \centering
    \begin{subfigure}{0.28\textwidth}
        \centering
        \includegraphics[width=\textwidth]{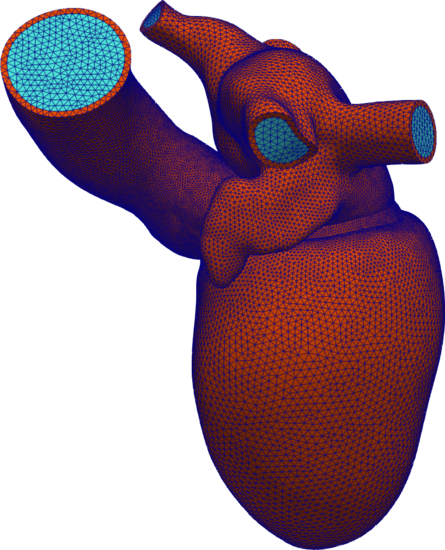}
        \caption{}
        \label{fig:mesh-solid}
    \end{subfigure}
    \begin{subfigure}{0.28\textwidth}
        \centering
        \includegraphics[width=\textwidth]{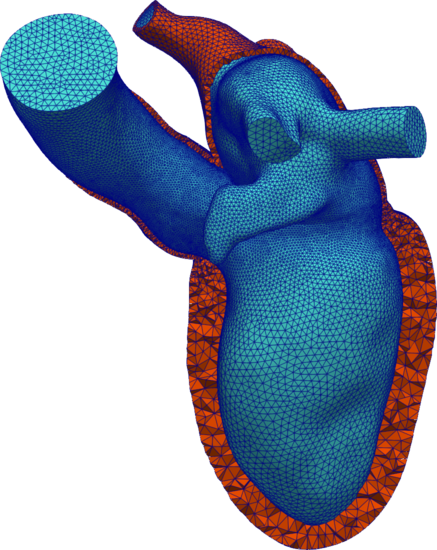}
        \caption{}
        \label{fig:mesh-fluid}
    \end{subfigure}
    \hspace{2em}
    \begin{subfigure}{0.225\textwidth}
        \centering
        \includegraphics[width=\textwidth]{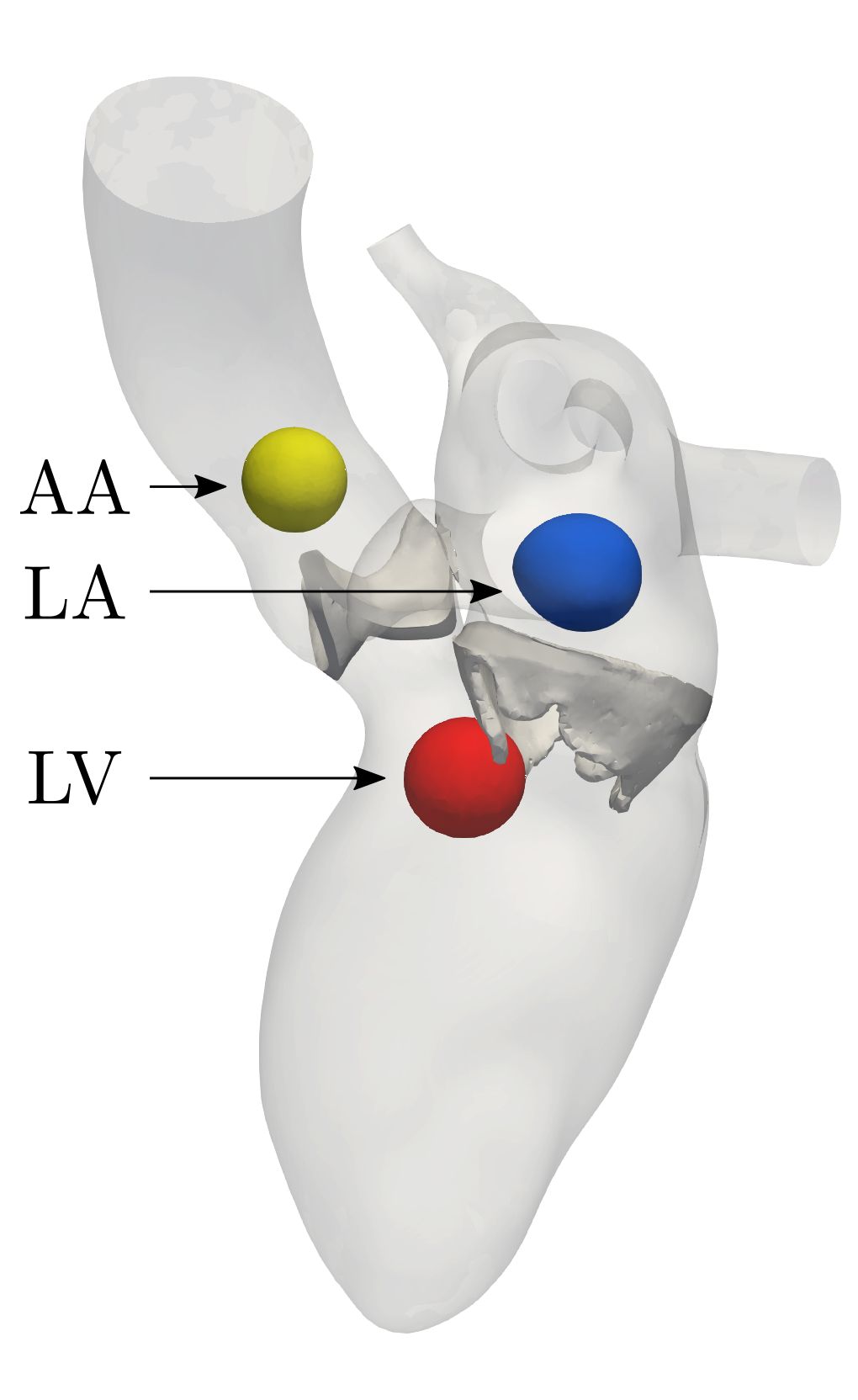}
        \caption{}
        \label{fig:control-volumes}
    \end{subfigure}

    \caption{(a, b) Computational mesh used for the solid (red) and the fluid (blue) domains. (c) Spherical control volumes used in estimating the average pressures within each chamber. Average pressure within control volumes is also used to trigger the opening and closing of valves.}
\end{figure}

\begin{table}
    \centering
    \begin{tabular}{c S[table-format=2.2] S[table-format=2.2] S[table-format=2.2] S[table-format=6.0] S[table-format=6.0]}
        \toprule
        & $h_{\min}$ [\si{\milli\metre}] & $h_\text{avg}$ [\si{\milli\metre}] & $h_{\max}$ [\si{\milli\metre}] & \textbf{\# elem.} & \textbf{\# nodes} \\
        \midrule
        \textbf{fluid} & 0.51 & 1.6 & 4.7 & 790533 & 137504 \\
        \textbf{solid} & 0.59 & 2.1 & 5.7 & 234132 & 63834 \\
        \bottomrule
    \end{tabular}
    \caption{Minimum, average and maximum mesh element diameter, number of elements and number of nodes used for spatial discretization.}
    \label{tab:mesh}
\end{table}

We introduce a tetrahedral mesh in both domains $\hat\Omega\solid$ and $\hat\Omega\fluid$, represented in \cref{fig:mesh-solid,fig:mesh-fluid}. At their interface $\hat\Sigma$, the fluid and solid meshes are conforming. The mesh $\hat\Omega\fluid$ is updated following the displacement $\mathbf d\ale$, resulting in a tetrahedral mesh over the domain $\Omega\fluid$. The spatial discretization is finer in the region near the valves, to allow better capturing their presence. We report details about the mesh size and number of elements in \cref{tab:mesh}.

The time-discrete ionic model \eqref{eq:ionic-time-discrete} is solved at every vertex in the solid mesh. Then, ionic current $I_\mathrm{ion}(v, \mathbf w, \mathbf z)$ is evaluated at every mesh vertex, and interpolated on quadrature nodes of the mesh. This approach is referred to as \textit{ionic current interpolation} (ICI) in the literature  \cite{krishnamoorthi2013numerical,pathmanathan2011significant}.

We discretize in space the electrophysiology model \eqref{eq:monodomain-time-discrete} with piecewise quadratic finite elements \cite{hughes2012finite,quarteroni2017numerical}. Indeed, quadratic finite elements have been shown to provide improved accuracy with a lower number of degrees of freedom, with respect to linear elements, for cardiac electrophysiology \cite{africa2022matrix}. Since the ionic current term is known from the solution of the ionic model, the resulting problem is linear and symmetric. We solve it with the conjugate gradient (CG) method \cite{quarteroni2017numerical,saad2003iterative}, preconditioned with an algebraic multigrid (AMG) preconditioner \cite{xu2017algebraic}. We point out that, in principle, while we do not consider it in this work, a staggered scheme may allow to use a finer spatial discretization for the electrophysiology problem, e.g. by using intergrid transfer operators as presented in \cite{salvador2020intergrid}.

The system \eqref{eq:activation-time-discrete} is solved at every node in the computational mesh. The regularization problem \eqref{eq:sl} is solved by means of linear finite elements, using a regularization radius $\delta_{SL}$ proportional to the mesh size. The resulting linear system is solved with the CG method, using AMG preconditioning.

The ALE lifting problem \eqref{eq:lifting-time-discrete} is discretized in space with piecewise linear finite elements. The resulting problem is linearized with Newton's method and then solved with the GMRES method \cite{quarteroni2017numerical,saad2003iterative}, preconditioned with AMG.

The FSI problem \eqref{eq:fsi-time-discrete} is discretized monolithically, using condensation of interface variables \cite{bucelli2022partitioned,gerbi2018numerical,nordsletten2011fluid,zhang2001analysis}. We found that the monolithic approach is computationally more efficient and robust with respect to partitioned approaches based on subiterations between the fluid and the solid models, as detailed in \cite{bucelli2022partitioned}. Piecewise linear finite elements are used for the discretization of fluid velocity, fluid pressure and solild displacement. The use of VMS-LES model for the fluid equations yields a stable numerical solution even though linear finite element spaces for pressure and velocity do not fulfill the inf-sup condition \cite{quarteroni2017numerical}. The VMS-LES model also provides stabilization for the advection-dominated regime \cite{bazilevs2007variational}. The resulting problem is linearized with Newton's method, and solved with the GMRES method, using a block lower-triangular preconditioner described in \cite{bucelli2022partitioned}.

\section{Numerical results}
\label{sec:results}

\begin{figure}
    \centering

    \begin{subfigure}{0.574\textwidth}
        \includegraphics[width=\textwidth]{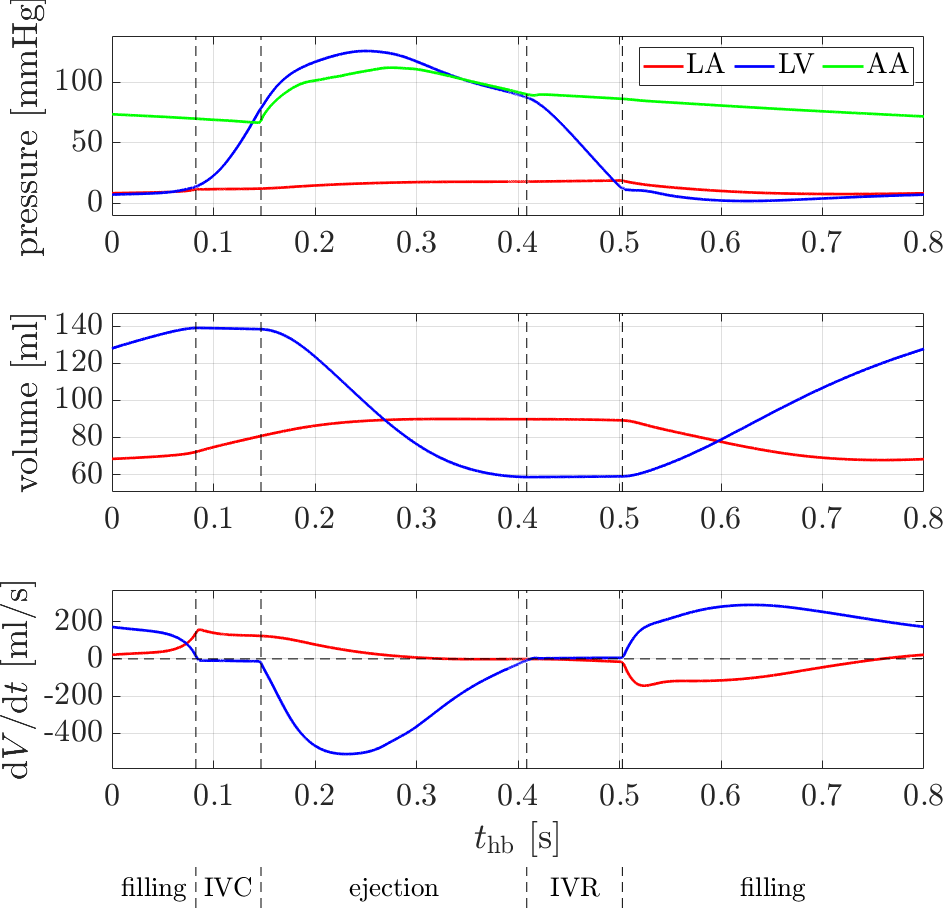}
        \caption{}
        \label{fig:pressure-volume-flowrate}
    \end{subfigure}
    \begin{subfigure}{0.39\textwidth}
        \includegraphics[width=\textwidth]{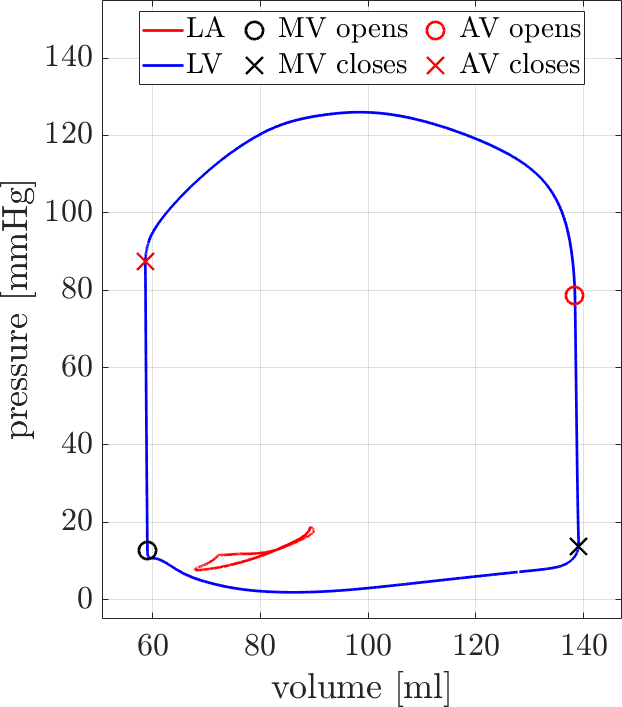}
        \caption{}
        \label{fig:pv-loop}
    \end{subfigure}
    \caption{(a) Evolution, over the third simulated heartbeat, of the pressure $p$, volume $V$ and volume derivative $\text{d}V/\text{dt}$ for the left heart compartments. Vertical dashed lines separate heartbeat phases, labelled on the bottom. (b) Pressure-volume loops for the left ventricle and atrium, with marks indicating opening and closing of valves.}
    \label{fig:pv-plots}
\end{figure}

We run numerical simulations on a left heart model provided by Zygote Media Group \cite{zygote}, representing the heart of an average 21-year-old male. Preprocessing of the geometry was done using the algorithms presented in \cite{fedele2021polygonal} as implemented in the software library \texttt{VMTK} \cite{vmtk}. Solvers for the individual core models and for the coupled model were implemented in the C++ library \lifex \cite{africa2022lifex,africa2022lifexcore,lifex}, based on the finite element core \dealii \cite{arndt2020dealii9.2,arndt2020dealii,dealii}. We run the simulations using 192 parallel processes on the GALILEO100 supercomputer\footnote{See \url{https://wiki.u-gov.it/confluence/display/SCAIUS/UG3.3\%3A+GALILEO100+UserGuide} for technical specifications.} at the CINECA high-performance computing center (Italy).

We simulate three heartbeats, setting $T = 3T_\text{hb} = \SI{2.4}{\second}$, and consider only the last heartbeat in our analysis (starting at $t_0 = 2T_\text{hb} = \SI{1.6}{\second}$), to reduce the effect of initial conditions. We shall denote with $t_\text{hb} = t - t_0$ times relative to the third heartbeat.

For each heartbeat, the simulation takes approximately \SI{21}{\hour} of wall time, of which \SI{93}{\percent} is spent assembling and solving the FSI problem, \SI{6}{\percent} assembling and solving the ALE lifting problem and \SI{2}{\percent} assembling and solving the electrophysiology problem.

We found that the computational time associated to the ALE lifting is considerably higher using the non-linear problem \eqref{eq:lifting} than using simpler harmonic or linear elastic lifting operators \cite{stein2003mesh} (for which the wall time would be less than \SI{1}{\percent} of the total). However, the latter frequently lead to mesh element inversion and solver breakdown, whereas the non-linear operator we employ has proven to be much more robust with respect to deformations.

For post-processing, we compute the average pressure of each compartment (LA, LV, AA) by averaging the pressure inside a spherical control volume within that compartment (see \cref{fig:control-volumes}). In \cref{fig:pv-plots}, we report the average pressure, volume and volume derivative over time of the compartments. \Cref{tab:indicators} collects some quantitative indicators for the heart function, comparing the values obtained by our simulation with data from the medical literature. Our numerical results are consistent with the clinically measured ranges. In the following sections, we provide details on the simulation results for each of the four heartbeat phases (identified in the plots of \cref{fig:pv-plots}).

\begin{table}
    \centering
    \begin{tabular}{l r S[table-format=3.1] S c l}
        \toprule
        \multicolumn{2}{c}{\textbf{Indicator}} & \textbf{Simulation} & \multicolumn{2}{c}{\textbf{Normal values}} & \textbf{Description} \\
        \midrule
        EDV                  & [\si{\milli\litre}]            & 139.1 & \numrange{126}{208} & \cite{maceira2006normalized} & left-ventricular end-diastolic volume \\
        ESV                  & [\si{\milli\litre}]            & 58.6  & \numrange{35}{80}   & \cite{maceira2006normalized} & left-ventricular end-systolic volume \\
        SV                   & [\si{\milli\litre}]            & 80.4  & \numrange{81}{137}  & \cite{maceira2006normalized} & left-ventricular stroke volume \\
        EF                   & [\si{\percent}]                & 57.8  & \numrange{49}{73}   & \cite{clay2006normal} & left-ventricular ejection fraction \\
        $p\lv_{\max}$        & [\si{\mmhg}]                   & 126.0   & 119(13) & \cite{sugimoto2017echocardiographic} & left-ventricular peak systolic pressure \\
        $Q^\text{AV}_{\max}$ & [\si{\milli\litre\per\second}] & 510.0   & 427(129) & \cite{hammermeister1974rate} & peak-systolic aortic flow rate \\
        $T_\text{IVC}$       & [\si{\milli\second}]           & 64.2  & \numrange{51}{90}   & \cite{fabian1972duration} & isovolumetric contraction time \\
        $T_\text{ej}$        & [\si{\milli\second}]           & 261.0   & \numrange{230}{334} & \cite{fabian1972duration} & ejection time \\
        $T_\text{IVR}$       & [\si{\milli\second}]           & 94.2  & \numrange{50}{140}  & \cite{benchimol1967study} & isovolumetric relaxation time \\
        $T_\text{fil}$       & [\si{\milli\second}]           & 379.0   & \numrange{280}{472} & \cite{little1990clinical} & diastolic filling time \\
        LFS                  & [\si{\percent}]                & 17.8  & \numrange{13}{21}   & \cite{emilsson2006mitral} & longitudinal fractional shortening \\
        \bottomrule
    \end{tabular}
    \caption{Values of physiological indicators computed from simulation results, and associated normal values from the medical literature. We report either normal ranges or mean $\pm$ standard deviation.}
    \label{tab:indicators}
\end{table}

\subsection{Isovolumetric contraction}

\begin{figure}
    \centering

    \begin{subfigure}{0.3548\figwidth}
        \centering
        \includegraphics[width=\textwidth]{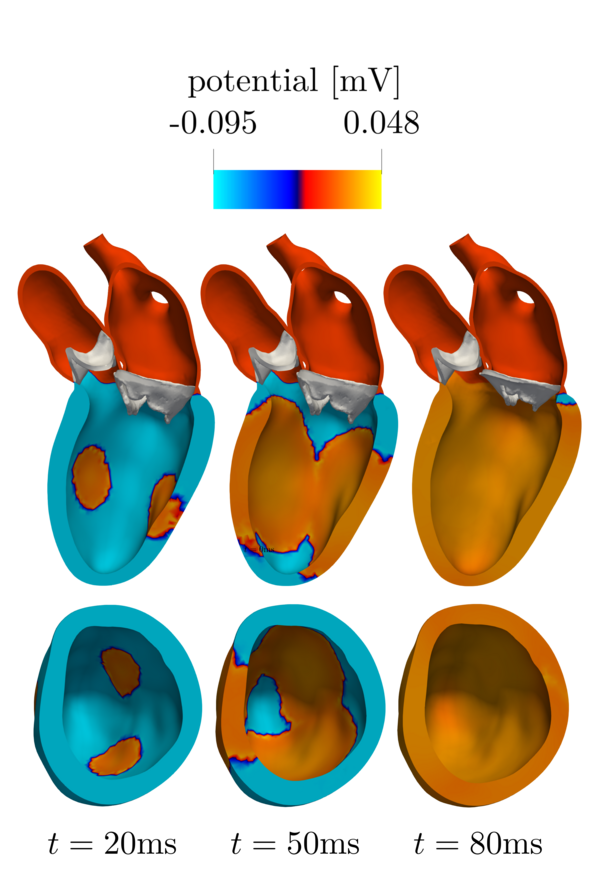}
        \caption{}
        \label{fig:ivc-potential}
    \end{subfigure}
    \begin{subfigure}{0.2448\figwidth}
        \centering
        \includegraphics[width=\textwidth]{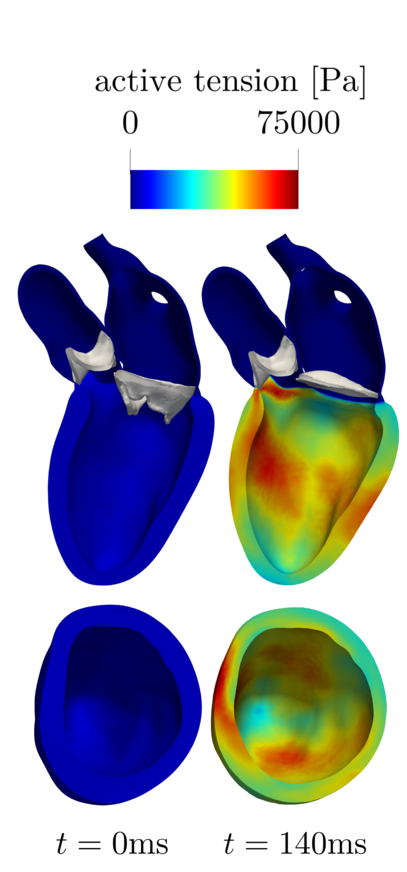}
        \caption{}
        \label{fig:ivc-active-stress}
    \end{subfigure}
    \begin{subfigure}{0.4003\figwidth}
        \centering
        \includegraphics[width=\textwidth]{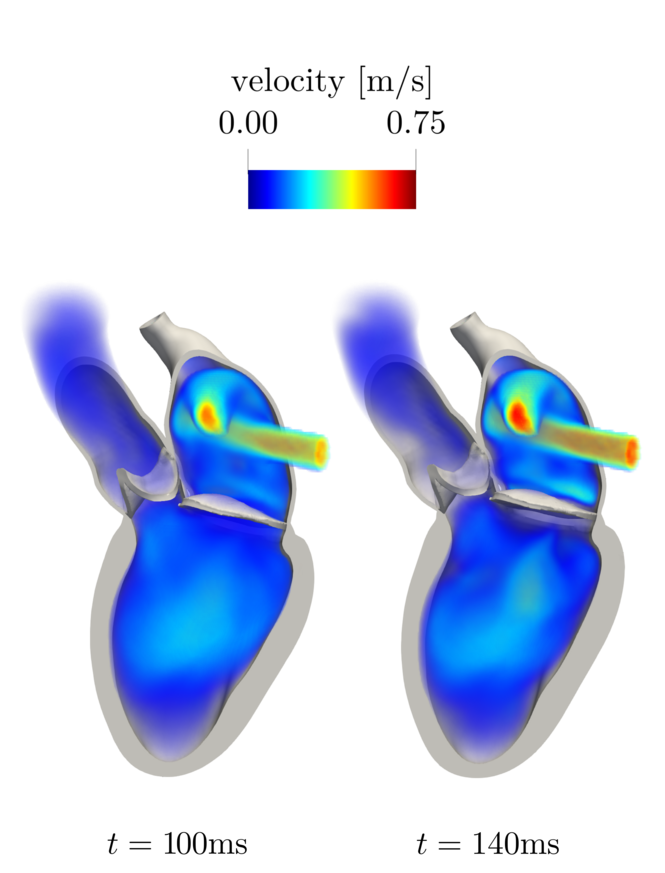}
        \caption{}
        \label{fig:ivc-velocity}
    \end{subfigure}

    \begin{subfigure}{0.3548\figwidth}
        \centering
        \includegraphics[width=\textwidth]{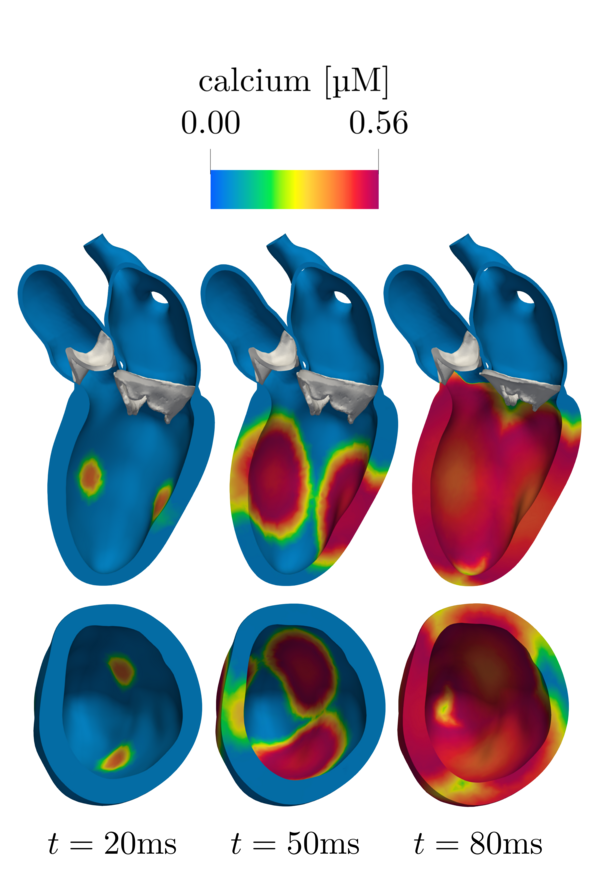}
        \caption{}
        \label{fig:ivc-calcium}
    \end{subfigure}
    \begin{subfigure}{0.2448\figwidth}
        \centering
        \includegraphics[width=\textwidth]{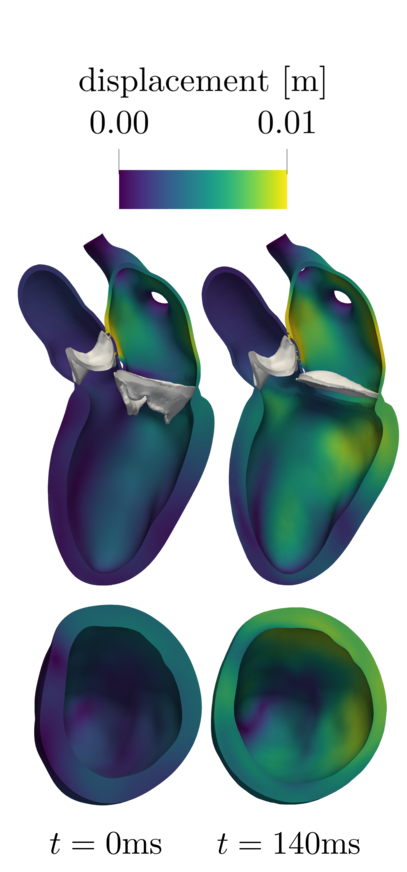}
        \caption{}
        \label{fig:ivc-displacement}
    \end{subfigure}
    \begin{subfigure}{0.4003\figwidth}
        \centering
        \includegraphics[width=\textwidth]{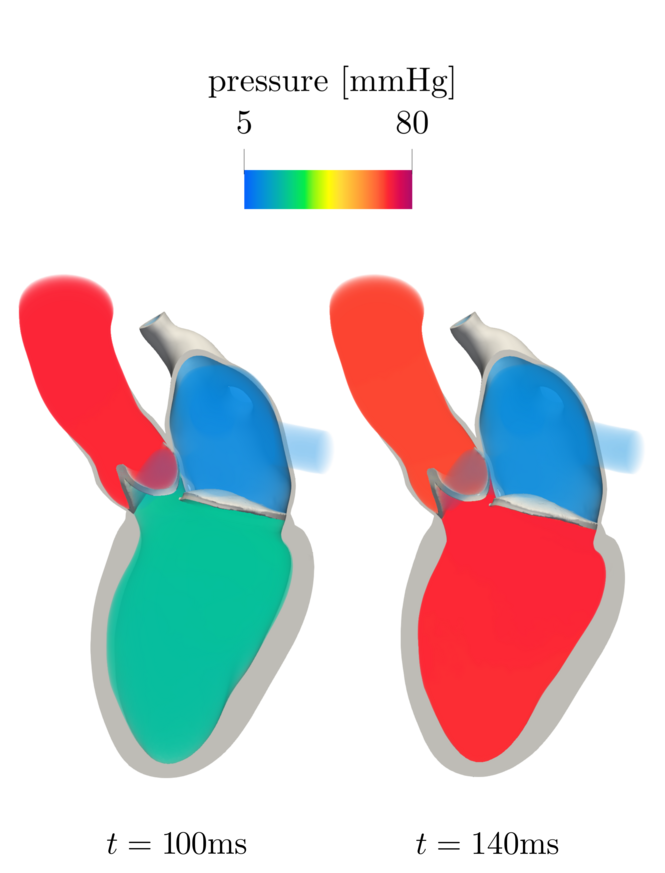}
        \caption{}
        \label{fig:ivc-pressure}
    \end{subfigure}

    \caption{Solution snapshots during the isovolumetric contraction phase, as seen through a long-axis and a short-axis section: (a) transmembrane potential $v$; (b) active tension $T_\text{act}$; (c) volume rendering of the fluid velocity magnitude $|\mathbf u|$; (d) intracellular calcium concentration $[\text{Ca}^{2+}]_\text{i}$; (e) displacement magnitude $|\mathbf d|$; (f) volume rendering of the fluid pressure $p$.}
    \label{fig:ivc}
\end{figure}

The heartbeat starts from the end of the diastolic filling phase, with the MV open and the AV closed (see \cref{fig:pv-plots}). The electrical stimulus is applied to the activation points of the left ventricle (see \cref{sec:electrophysiology}), leading to it being completely activated within $t_\text{hb} = \SI{100}{\milli\second}$ (\cref{fig:ivc-potential}). This leads to the increase of intracellular calcium concentration (\cref{fig:ivc-calcium}), which subsequently determines the generation of contractile force within the ventricular wall (\cref{fig:ivc-active-stress}). Intraventricular pressure rises, triggering the closure of the MV. When the valve is closed ($t_\text{hb} = \SI{89}{\milli\second}$), the ventricle is at its \textit{end-diastolic volume} $\text{EDV} = \SI{139}{\milli\litre}$.

At this point, the \textit{isovolumetric contraction} (IVC) phase starts: both valves are closed, and the ventricular pressure increases rapidly (see \cref{fig:pv-plots} and \cref{fig:ivc-pressure}) while the volume is maintained approximately constant.
Our model does not capture exactly the conservation of ventricular blood volume during this phase, due to the explicit discretization of the fluid domain displacement \eqref{eq:lifting-time-discrete} and to the use of the resistive model for valves, which allows for a little flow through the immersed surfaces. Nonetheless, during isovolumetric contraction the maximum volume variation equals \SI{0.7}{\milli\litre}, corresponding to \SI{0.5}{\percent} of the EDV. We deem this spurious variation to be acceptable, in accordance with similar spurious variations observed in the cardiac modeling literature \cite{viola2020fluid}.

During the IVC phase, the ventricle undergoes a small deformation, with its shape becoming slightly more spherical (see \cref{fig:ivc-active-stress,fig:ivc-displacement}) as described in \cite{hawthorne1966dynamic,katz2010physiology} .

The IVC phase lasts for $T_\text{IVC} = \SI{64.2}{\milli\second}$, consistently with physiological behavior \cite{fabian1972duration,feher2017quantitative}. When the pressure in the ventricle becomes larger than the pressure in the aorta, the opening of the AV is triggered and the ejection phase starts.

We remark that the possibility of including isovolumetric phases in a three-dimensional hemodynamics model is distinctive of FSI \cite{bucelli2022partitioned}, as those phases cannot be represented by electromechanics-driven CFD models \cite{zingaro2022geometric}, unless ad-hoc techniques are implemented \cite{this2020augmented,zingaro2022modeling}.

\subsection{Ejection}

\begin{figure}
    \centering

    \begin{subfigure}{0.2404\figwidth}
        \centering
        \includegraphics[width=\textwidth]{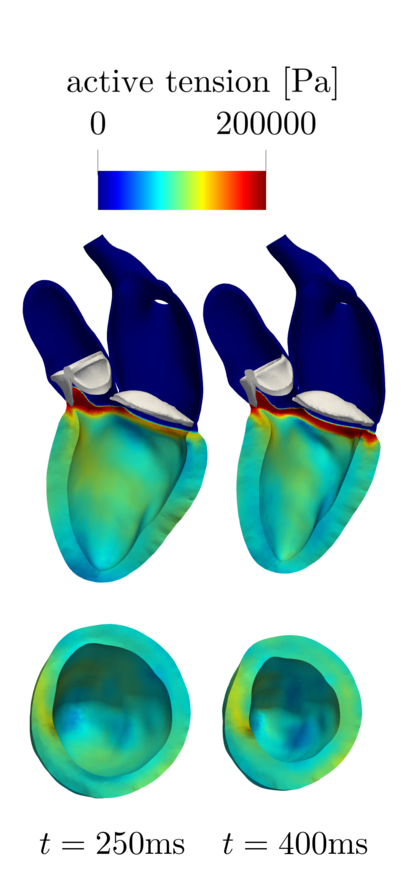}
        \caption{}
        \label{fig:ejection-active-stress}
    \end{subfigure}
    \begin{subfigure}{0.7596\figwidth}
        \centering
        \includegraphics[width=\textwidth]{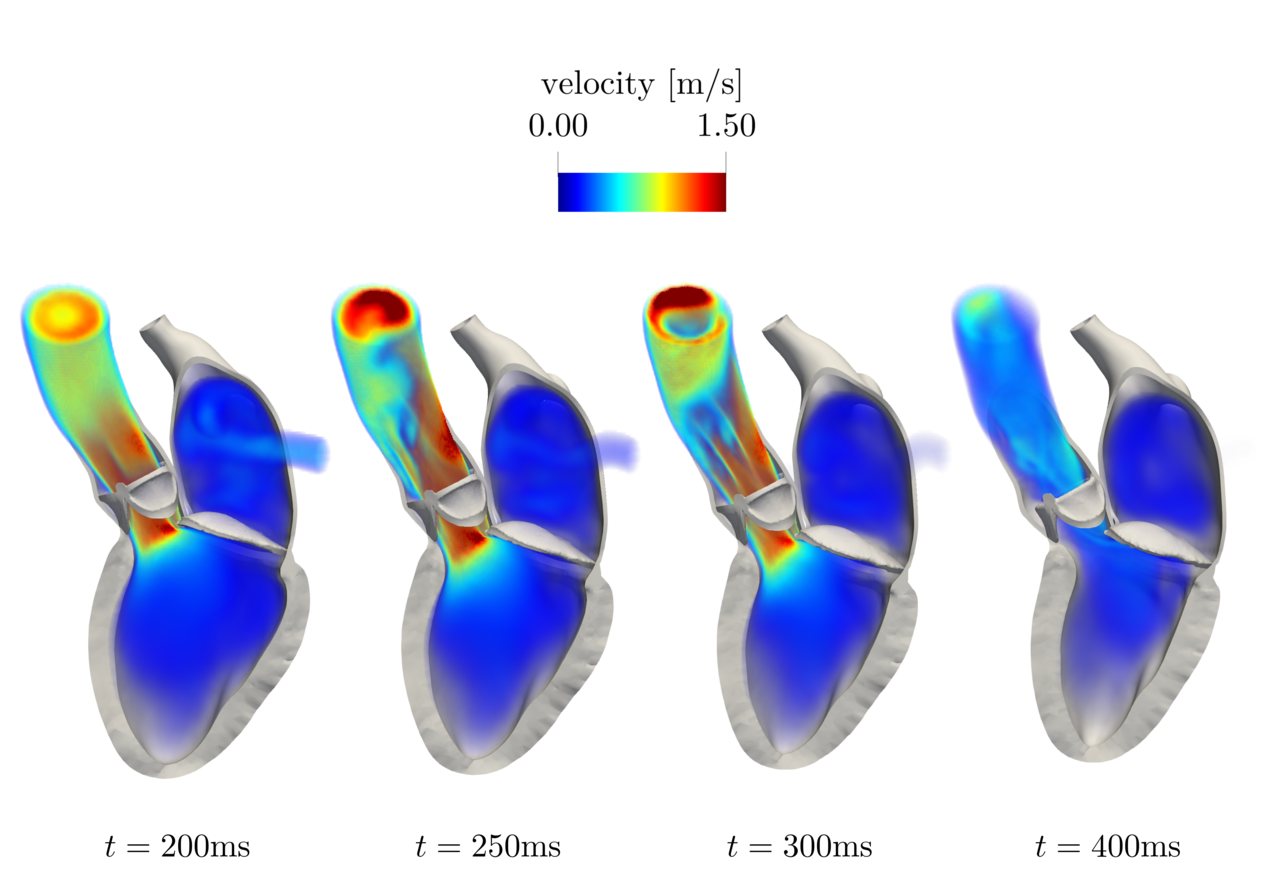}
        \caption{}
        \label{fig:ejection-velocity}
    \end{subfigure}

    \begin{subfigure}{0.2404\figwidth}
        \centering
        \includegraphics[width=\textwidth]{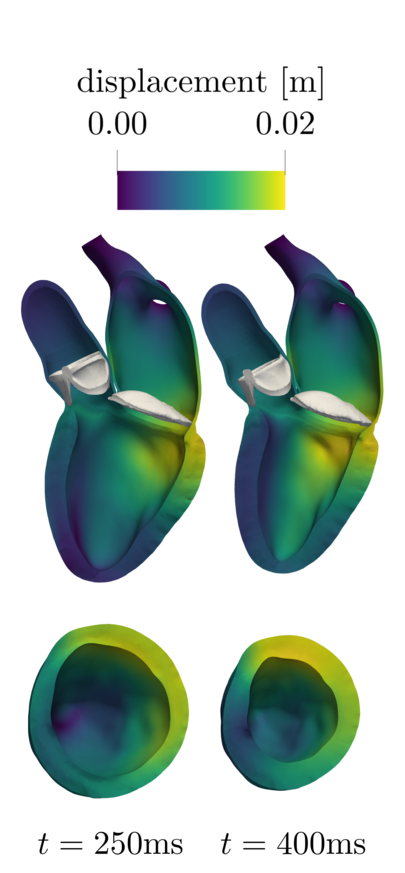}
        \caption{}
        \label{fig:ejection-displacement}
    \end{subfigure}
    \begin{subfigure}{0.7596\figwidth}
        \centering
        \includegraphics[width=\textwidth]{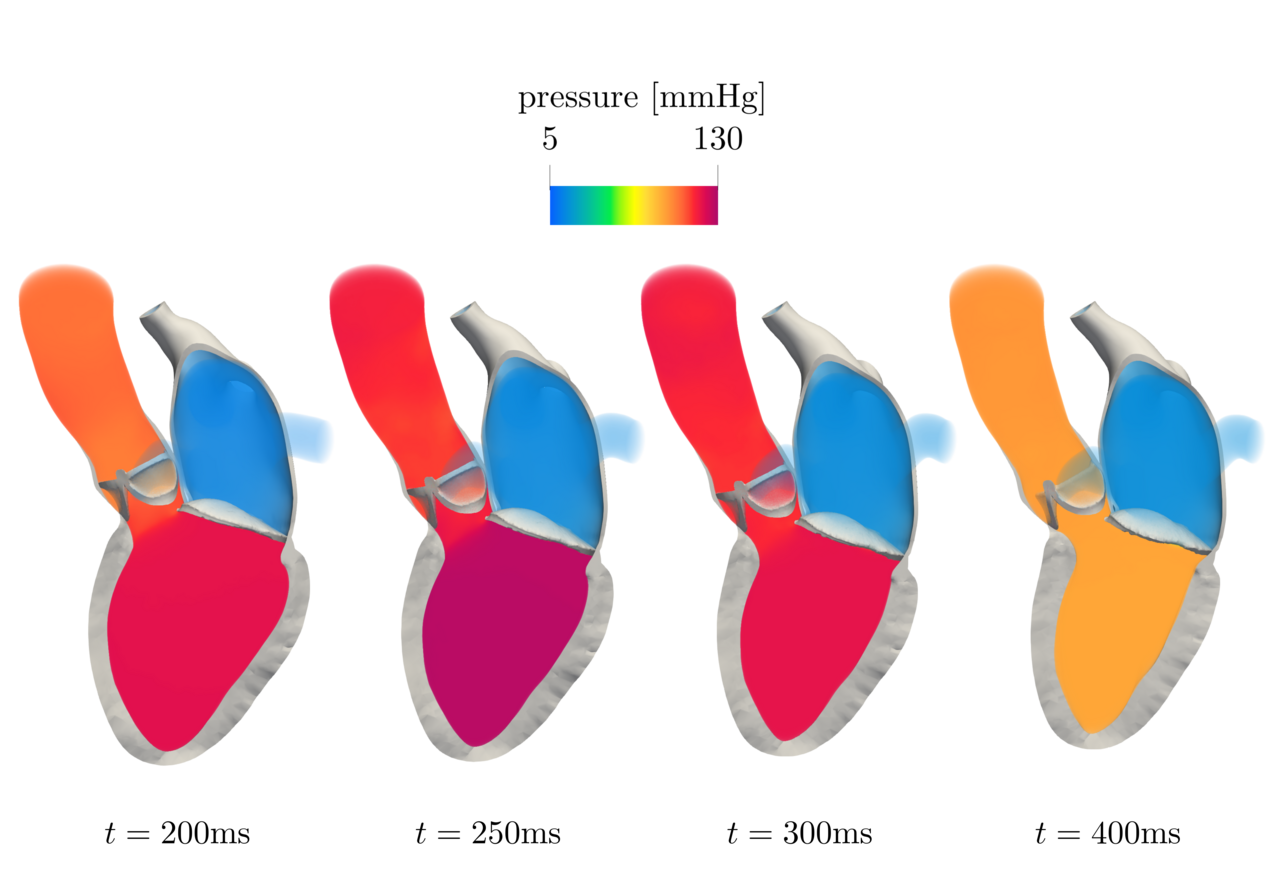}
        \caption{}
        \label{fig:ejection-pressure}
    \end{subfigure}

    \caption{Solution snapshots during the ejection phase, as seen through a long-axis and a short-axis section: (a) active tension $T_\text{act}$; (b) volume rendering of the fluid velocity magnitude $|\mathbf u|$; (c) displacement magnitude $|\mathbf d|$; (d) volume rendering of the fluid pressure $p$.}
    \label{fig:ejection}
\end{figure}

\begin{figure}
    \centering
    \includegraphics[width=0.4\textwidth]{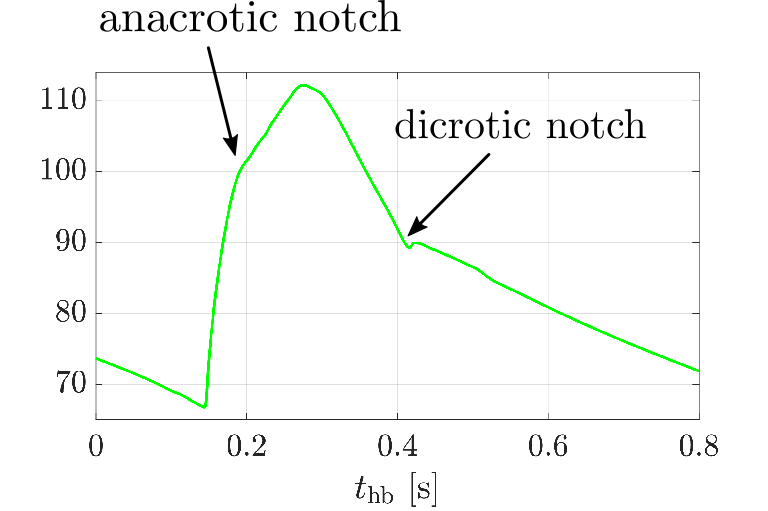}

    \caption{Average pressure in the ascending aorta over time. Notice the presence of the anacrotic notch during the systolic upstroke and of the dicrotic notch, or incisure, at the closing time of the AV.}
    \label{fig:aortic-pressure}
\end{figure}

\begin{figure}
    \centering

    \begin{subfigure}{0.35\textwidth}
        \includegraphics[width=\textwidth]{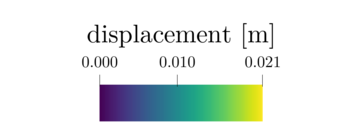}
    \end{subfigure}
    \vspace{0.5em}

    \begin{subfigure}{0.2\textwidth}
        \centering
        \includegraphics[width=\textwidth]{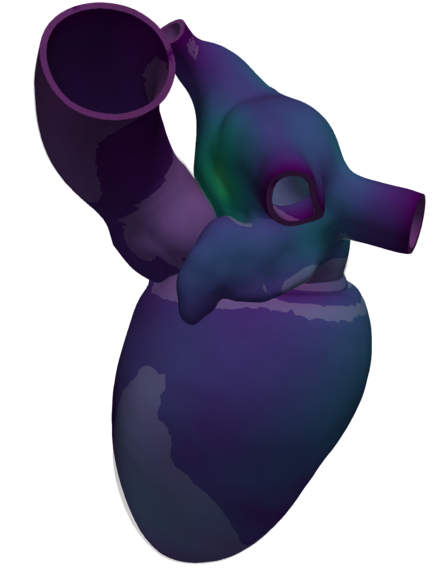}
         \caption{$t_\text{hb} = \SI{120}{\milli\second}$}
    \end{subfigure}
    \begin{subfigure}{0.2\textwidth}
        \centering
        \includegraphics[width=\textwidth]{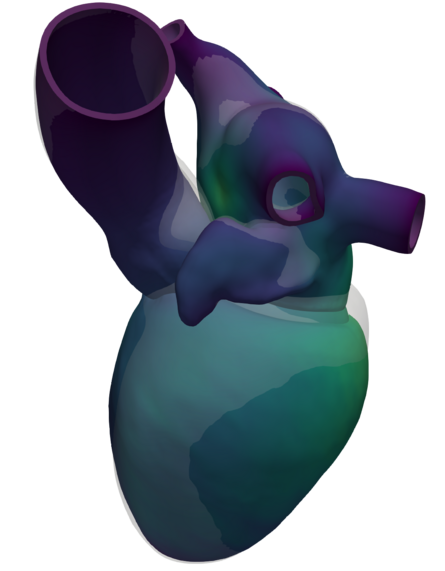}
         \caption{$t_\text{hb} = \SI{190}{\milli\second}$}
    \end{subfigure}
    \begin{subfigure}{0.2\textwidth}
        \centering
        \includegraphics[width=\textwidth]{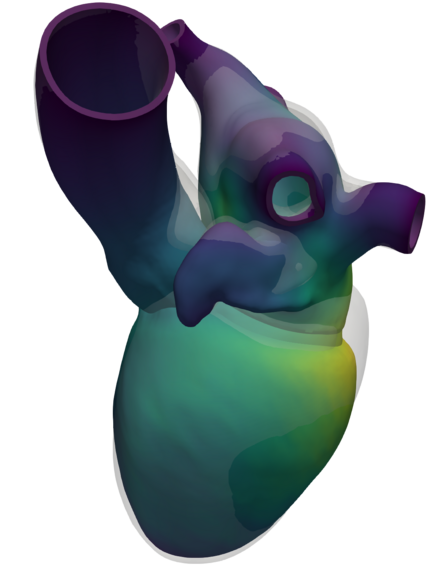}
         \caption{$t_\text{hb} = \SI{260}{\milli\second}$}
    \end{subfigure}
    \begin{subfigure}{0.2\textwidth}
        \centering
        \includegraphics[width=\textwidth]{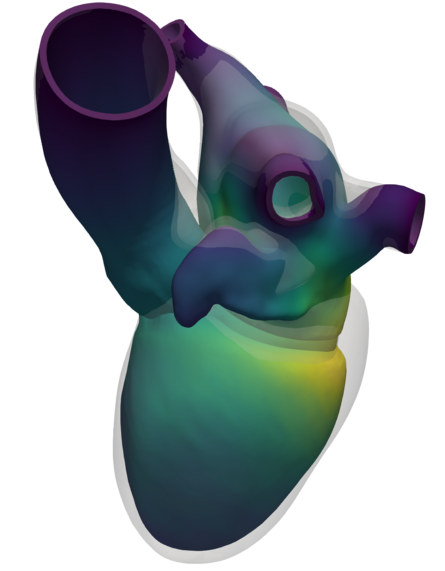}
         \caption{$t_\text{hb} = \SI{330}{\milli\second}$}
    \end{subfigure}
    \caption{Magnitude of the displacement $\mathbf d$ at four instants during the ejection phase. The initial configuration is overlaid in trasparency, and the domain is warped by $\mathbf d$. Notice how the ventricle becomes shorter during the ejection, mostly due to the shift of the atrioventricular plane towards the apex.}
    \label{fig:shortening}
\end{figure}

Blood is ejected from the ventricle into the aorta (\cref{fig:ejection-velocity}). The maximum flow rate through the AV orifice is attained at $t_\text{hb} = \SI{231}{\milli\second}$, and it equals $Q_{\max}^\text{AV} = \SI{510}{\milli\litre\per\second}$, consistently with physiology \cite{gulsin2017cardiovascular,hammermeister1974rate}. It corresponds to an average velocity magnitude of $|\mathbf u|_{\max}^\text{AV} = \SI{1.46}{\meter\per\second}$ on the AV section. The jet through the aortic valve is unsteady during all the ejection, consistently with the transitional nature of the flow (see \cref{fig:ejection-velocity}). At the end of the ejection phase, the left ventricle attains its \textit{end-systolic volume} $\text{ESV} = \SI{58.6}{\milli\litre}$. Given EDV and ESV, we compute the \textit{stroke volume} SV and \textit{ejection fraction} EF as
\[
\text{SV} = \text{EDV} - \text{ESV} = \SI{80.4}{\milli\litre} \qquad\qquad
\text{EF} = \frac{\text{SV}}{\text{EDV}} = \SI{57.8}{\percent}\;.
\]
Both quantities are within normal physiological ranges (see \cref{tab:indicators}) \cite{clay2006normal,feher2017quantitative,kumar2014robbins,maceira2006normalized,stanfield2016principles}.

During the ejection phase, the ventricular pressure increases until $t_\text{hb} = \SI{250}{\milli\second}$, reaching a peak value  of $p\lv_{\max} = \SI{126}{\mmhg}$ (see \cref{fig:pv-plots}), within the physiological ranges \cite{katz2010physiology,sugimoto2017echocardiographic}. After that, pressure starts decreasing until it falls below the aortic pressure, at which point the AV starts closing.

A similar evolution characterizes the pressure in the ascending aorta, reported in \cref{fig:aortic-pressure}: starting from an end-diastolic value of \SI{68}{\mmhg}, it reaches a peak of \SI{112}{\mmhg} and then declines until the next heartbeat. Although the absolute pressure values are smaller than normal ones \cite{mark1998atlas}, and there is a large pressure jump between the ventricle and the aorta (see \cref{sec:limitations}), the time profile of the aortic pressure is remarkably similar to the ones obtained from in-vivo measurements \cite{mark1998atlas,murgo1980aortic}. In particular, it features the \textit{anacrotic notch}, resulting from the interaction of the forward and reflected pressure waves \cite{murgo1980aortic}. This effect is captured thanks to the FSI modeling framework, which, contrary to standalone CFD models, allows to obtain traveling pressure waves. The aortic pressure also features the \textit{dicrotic notch} \cite{katz2010physiology,mark1998atlas,stanfield2016principles} in correspondence of the AV closure.

As the volume reduces, the ventricle becomes shorter and the atrioventricular plane shifts towards the ventricular apex (see \cref{fig:shortening}), as observed in healthy hearts \cite{dusch2014diastolic,katz2010physiology,levrero2020sensitivity}. We quantify this effect by computing the \textit{longitudinal fractional shortening} (LFS) \cite{levrero2020sensitivity}: denoting by $L_\text{ED}$ and $L_\text{ES}$ the apico-basal distances at the end of diastole and at the end of systole, we have
\[
\text{LFS} = \frac{L_\text{ED} - L_\text{ES}}{L_\text{ED}} = \SI{17.8}{\percent}\;,
\]
which matches measurements on healthy hearts \cite{emilsson2006mitral}. At the end of the systolic phase, the ventricular wall is approximately \SI{14}{\percent} thicker than at the end of diastole \cite{pandian1983heterogeneity}.

\begin{figure}
    \centering

    \includegraphics[width=0.3\figwidth]{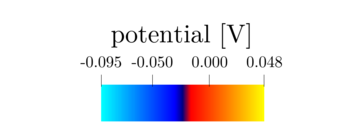}

    \hspace{1cm}

    \begin{subfigure}{0.2\figwidth}
        \centering
        \includegraphics[width=0.16\figwidth, trim={1in, 0in, 1in, 0in}, clip]{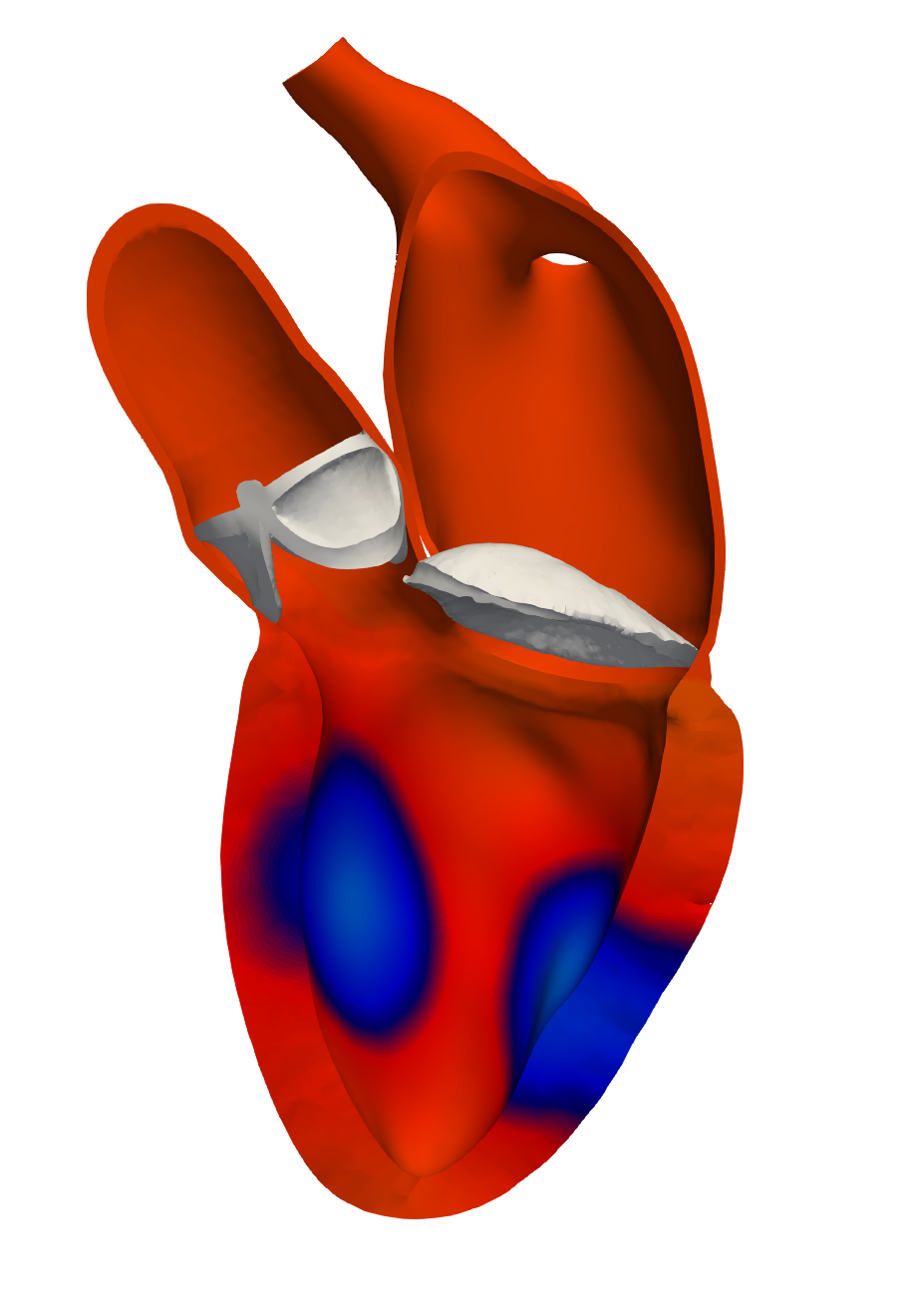}
        \caption{$t = \SI{300}{\milli\second}$}
    \end{subfigure}
    \begin{subfigure}{0.16\figwidth}
        \centering
        \includegraphics[width=0.16\figwidth, trim={1in, 0in, 1in, 0in}, clip]{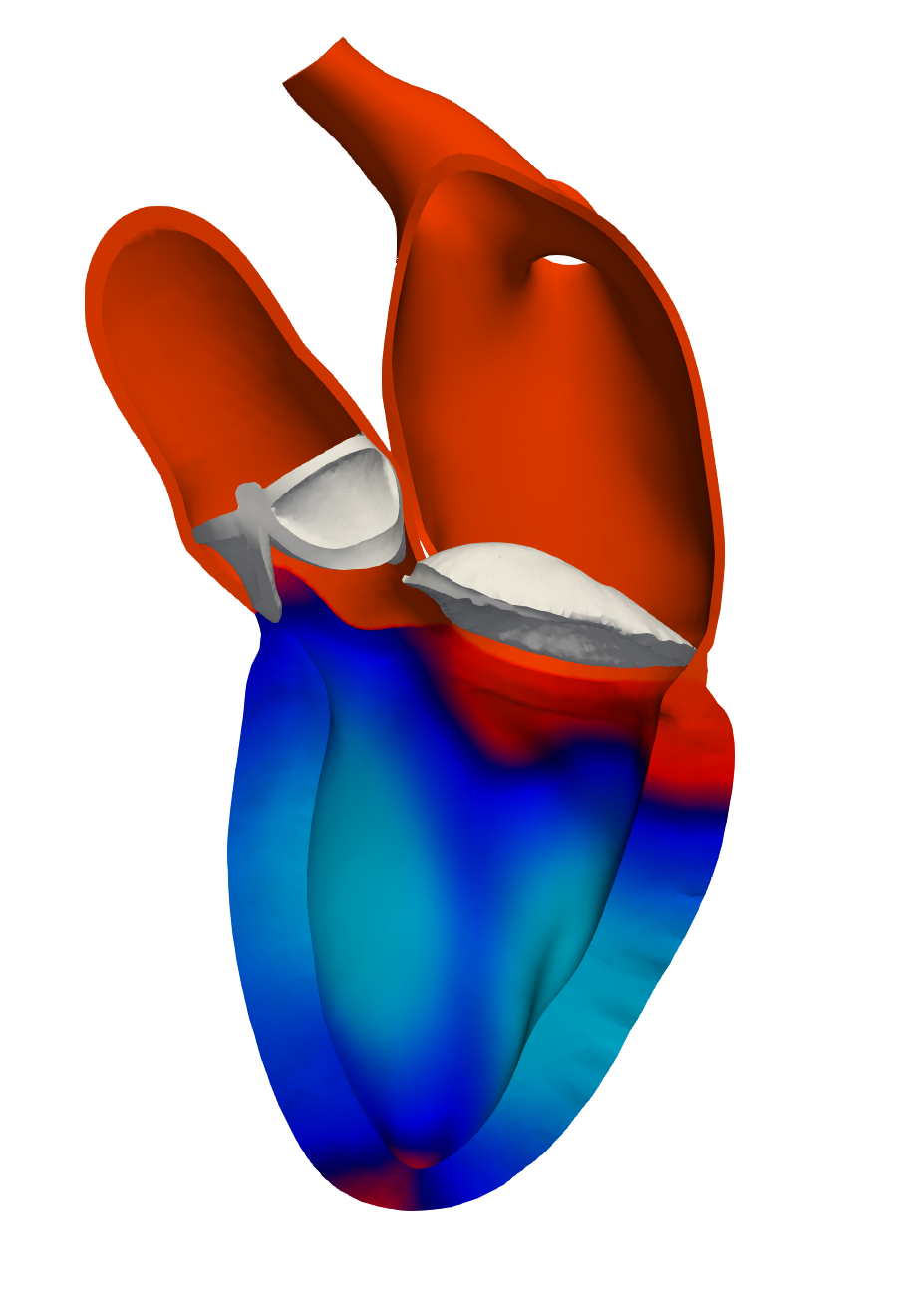}
        \caption{$t = \SI{330}{\milli\second}$}
    \end{subfigure}
    \begin{subfigure}{0.16\figwidth}
        \centering
        \includegraphics[width=0.16\figwidth, trim={1in, 0in, 1in, 0in}, clip]{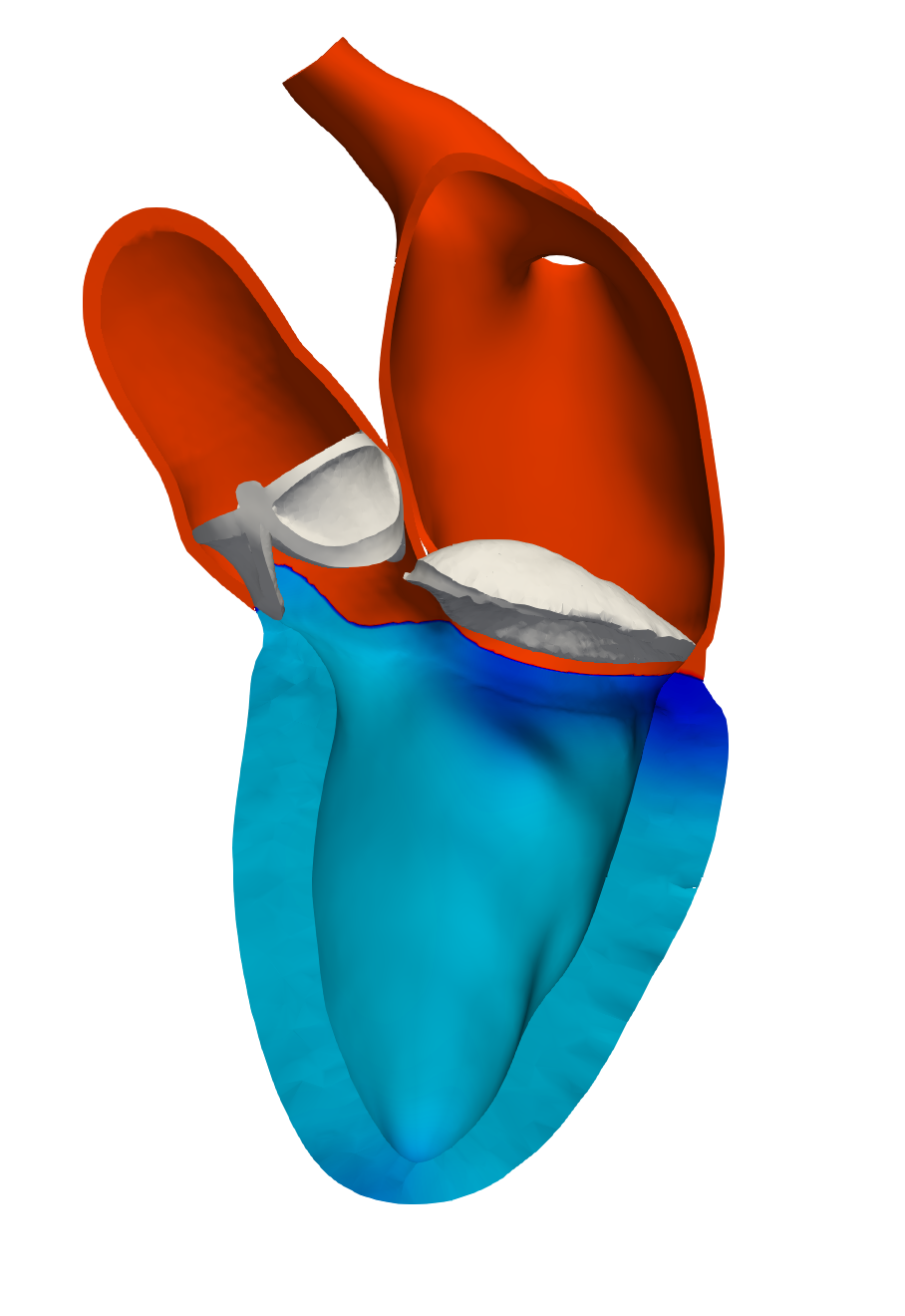}
        \caption{$t = \SI{360}{\milli\second}$}
    \end{subfigure}
    \begin{subfigure}{0.16\figwidth}
        \centering
        \includegraphics[width=0.16\figwidth, trim={1in, 0in, 1in, 0in}, clip]{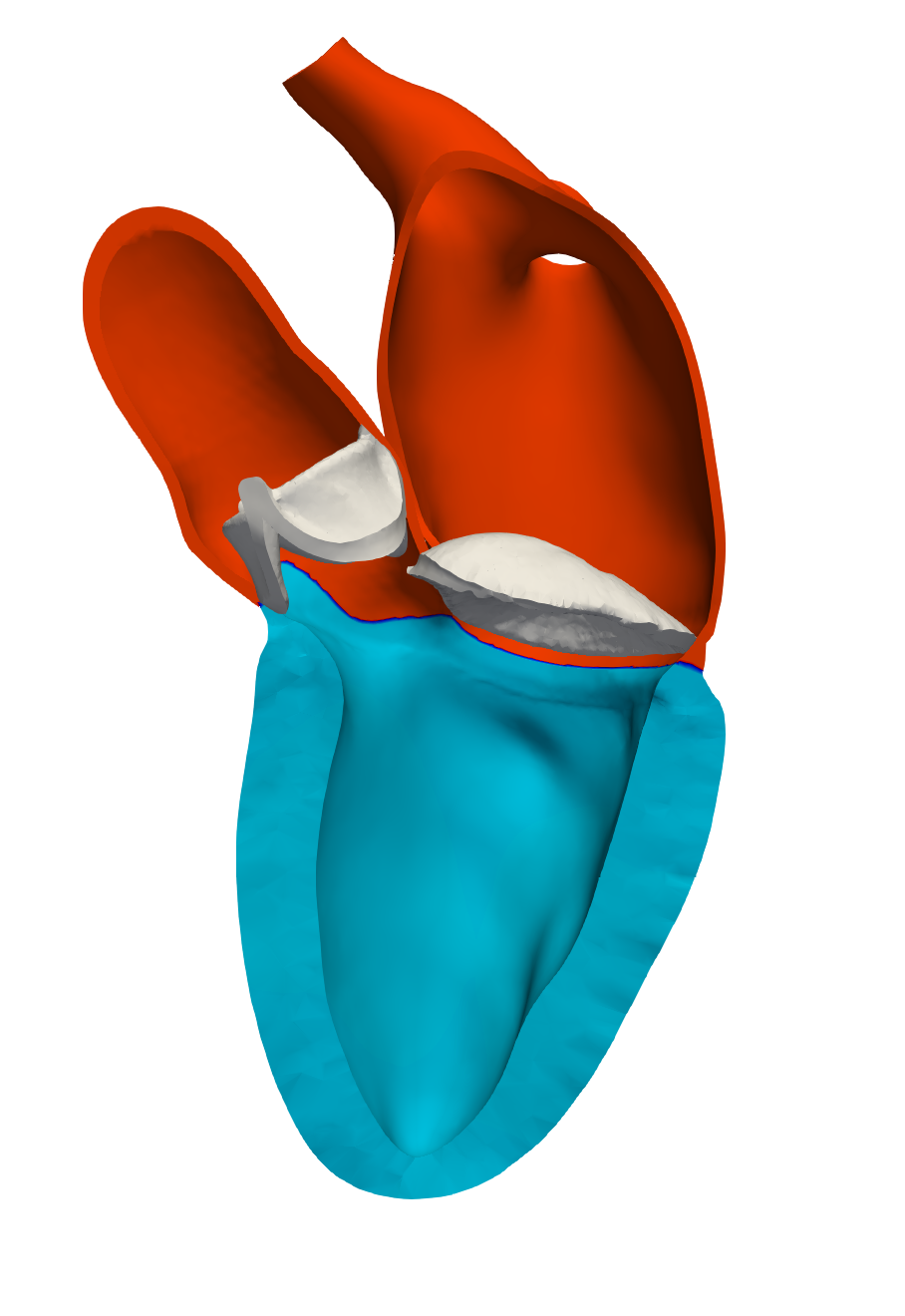}
        \caption{$t = \SI{400}{\milli\second}$}
    \end{subfigure}

    \caption{Repolarization of the left ventricle. The ventricular transmembrane potential gradually returns to its resting value.}
    \label{fig:repolarization}
\end{figure}

During the late systolic phase, towards the end of the ejection, the ventricle repolarizes (\cref{fig:repolarization}), with the transmembrane potential returning to its resting value.

Overall, the ejection phase lasts $T_\text{ej} = \SI{261}{\milli\second}$, and the whole systolic phase lasts $T_\text{sys} = \SI{326}{\milli\second}$, corresponding to \SI{40.7}{\percent} of the heartbeat.

\subsection{Isovolumetric relaxation}

\begin{figure}
    \centering

    \begin{subfigure}{0.2352\figwidth}
        \centering
        \includegraphics[width=\textwidth]{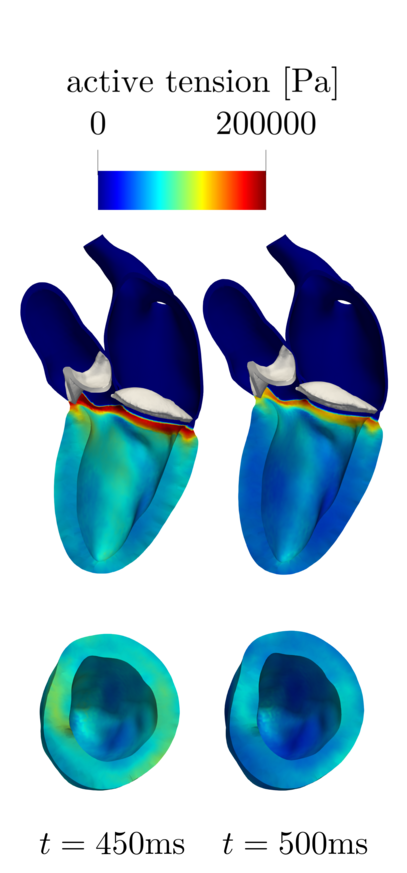}
        \caption{}
        \label{fig:ivr-active-stress}
    \end{subfigure}
    \begin{subfigure}{0.5648\figwidth}
        \centering
        \includegraphics[width=\textwidth]{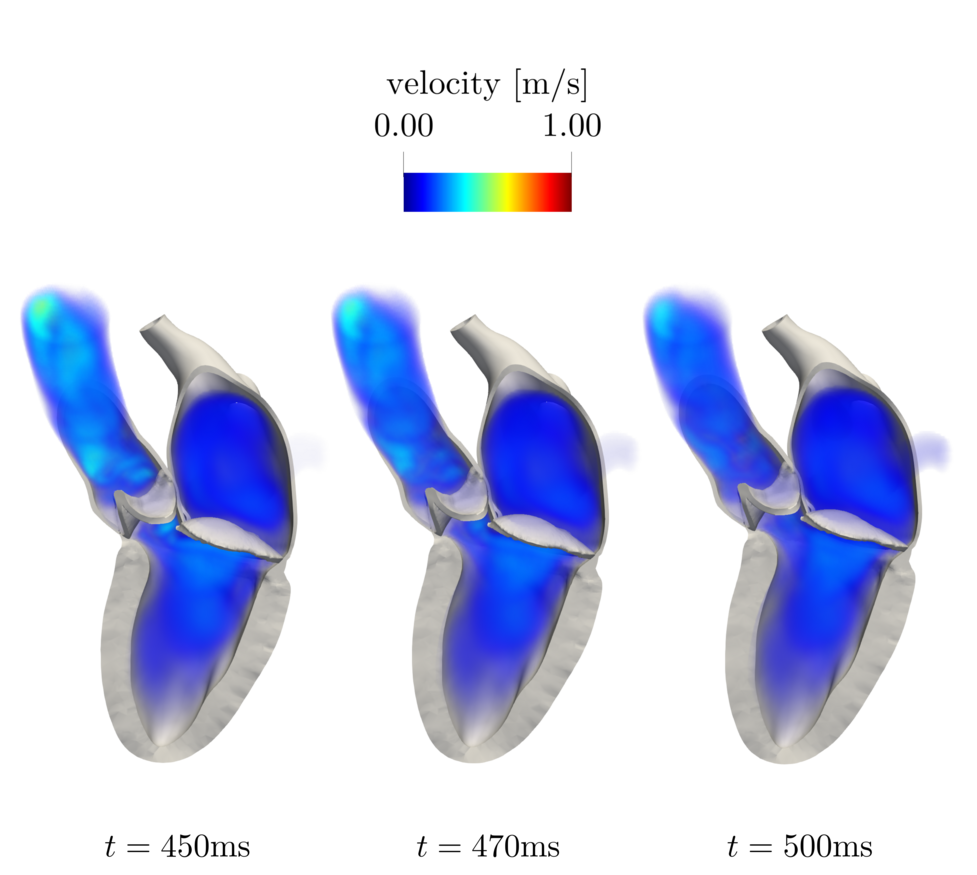}
        \caption{}
        \label{fig:ivr-velocity}
    \end{subfigure}

    \begin{subfigure}{0.2352\figwidth}
        \centering
        \includegraphics[width=\textwidth]{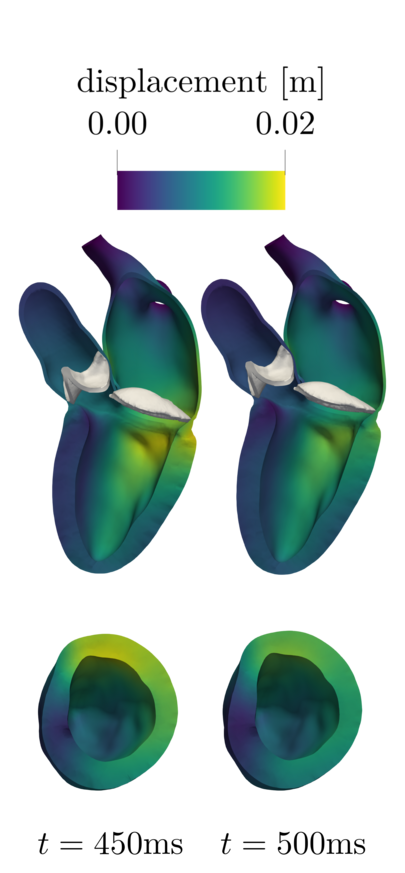}
        \caption{}
        \label{fig:ivr-displacement}
    \end{subfigure}
    \begin{subfigure}{0.5648\figwidth}
        \centering
        \includegraphics[width=\textwidth]{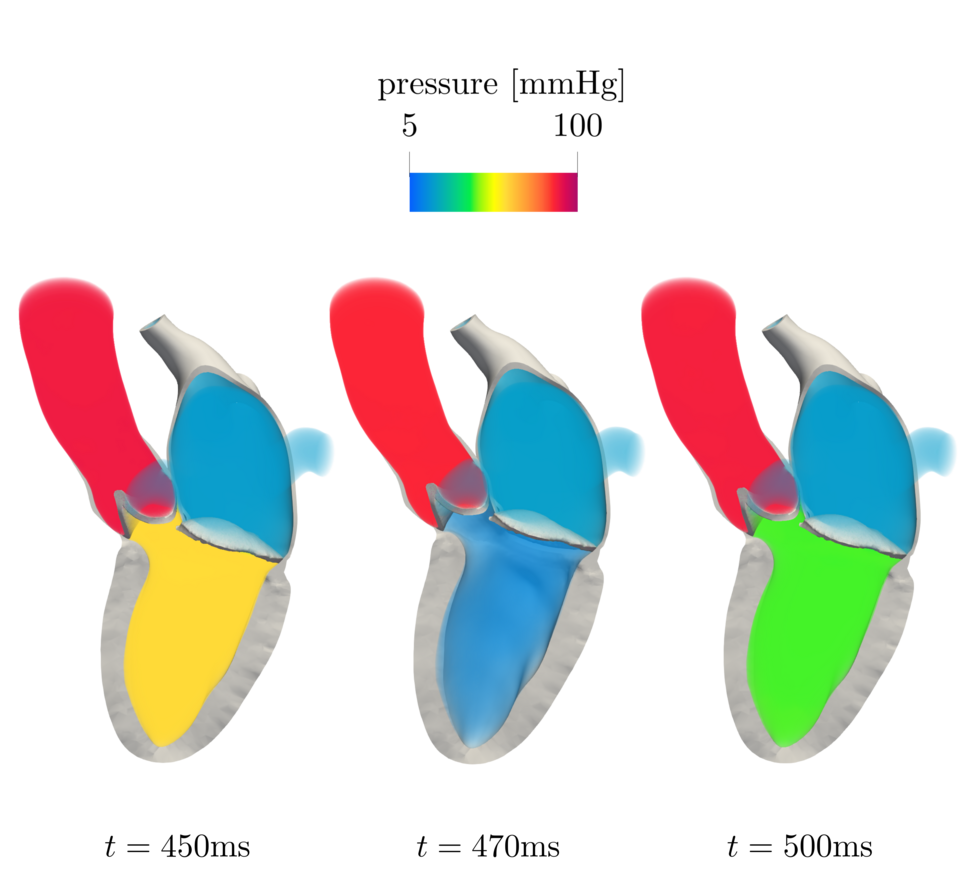}
        \caption{}
        \label{fig:ivr-pressure}
    \end{subfigure}

    \caption{Solution snapshots during the isovolumetric relaxation phase, as seen through a long-axis and a short-axis section: (a) active tension $T_\text{act}$; (b) volume rendering of the fluid velocity magnitude $|\mathbf u|$; (c) displacement magnitude $|\mathbf d|$; (d) volume rendering of the fluid pressure $p$.}
    \label{fig:ivr}
\end{figure}

Once the AV is fully closed, the \textit{isovolumetric relaxation} (IVR) phase starts (\cref{fig:ivr}). Ventricular pressure reduces as the ventricle relaxes at constant volume (\cref{fig:pv-plots}). This phase lasts for $T_\text{IVR} = \SI{94.2}{\milli\second}$, consistently with physiology \cite{benchimol1967study,little1990clinical}, and the MV starts opening as soon as the ventricular pressure becomes smaller than the atrial pressure.

As observed for the IVC phase, our model features a small spurious variation in volume during isovolumetric phases. In the case of the IVR phase, the variation amounts to \SI{0.4}{\milli\litre}, corresponding to \SI{0.6}{\percent} of the ESV. Also in this case, we deem the spurious variation to be acceptable.

\subsection{Filling}

\begin{figure}
    \centering

    \begin{subfigure}{0.2404\figwidth}
        \centering
        \includegraphics[width=\textwidth]{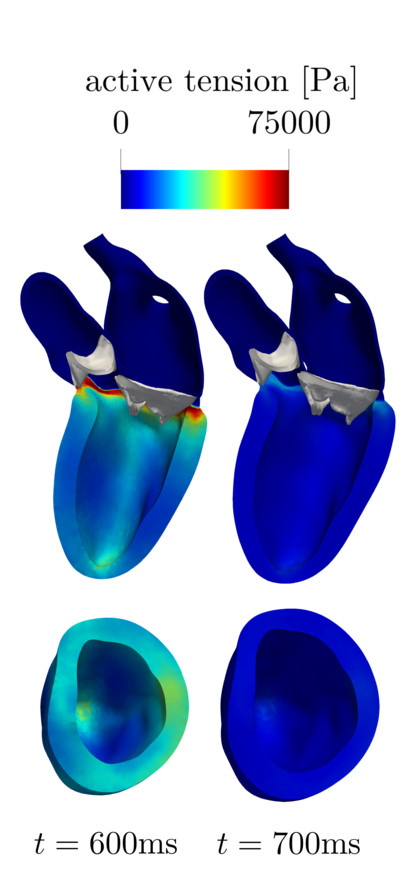}
        \caption{}
        \label{fig:filling-active-stress}
    \end{subfigure}
    \begin{subfigure}{0.7596\figwidth}
        \centering
        \includegraphics[width=\textwidth]{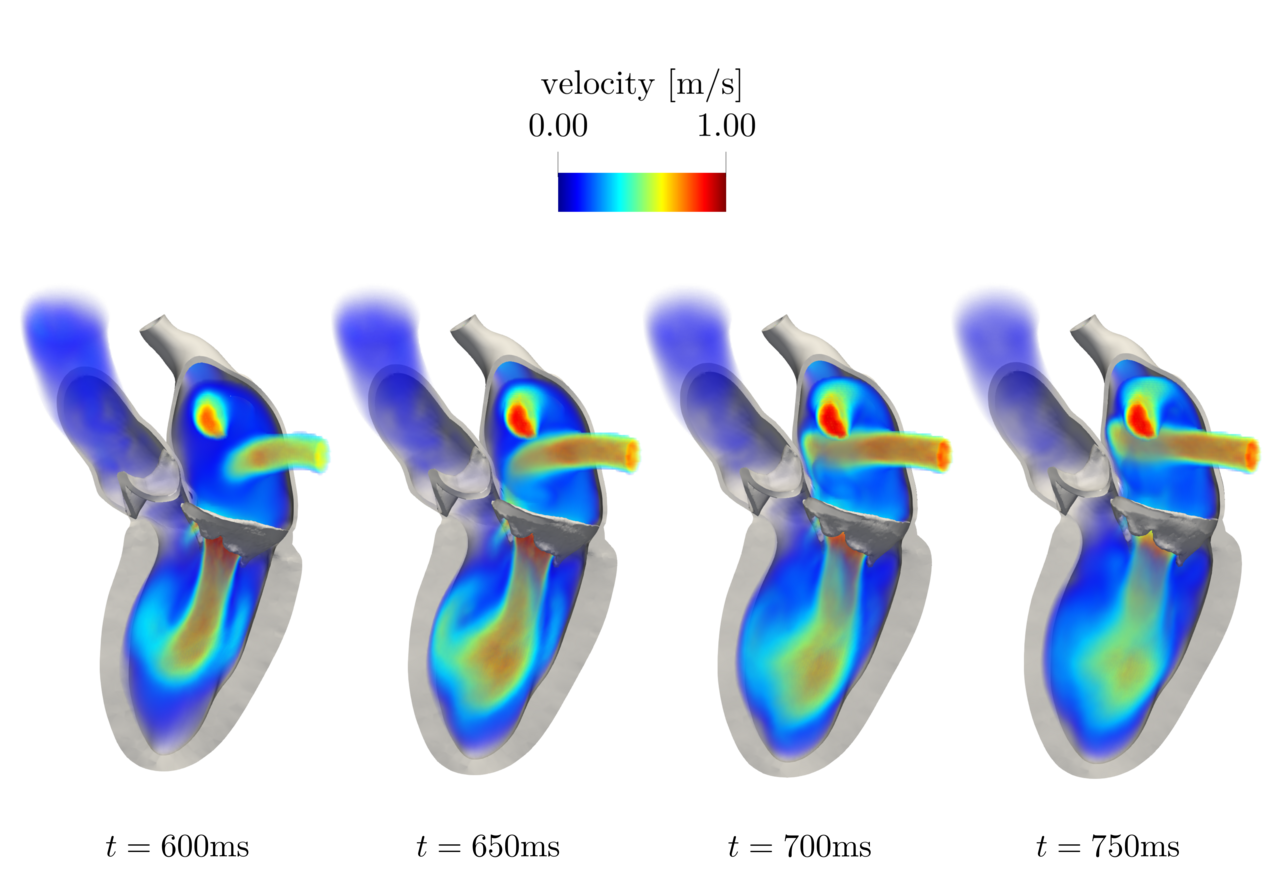}
        \caption{}
        \label{fig:filling-velocity}
    \end{subfigure}

    \begin{subfigure}{0.2404\figwidth}
        \centering
        \includegraphics[width=\textwidth]{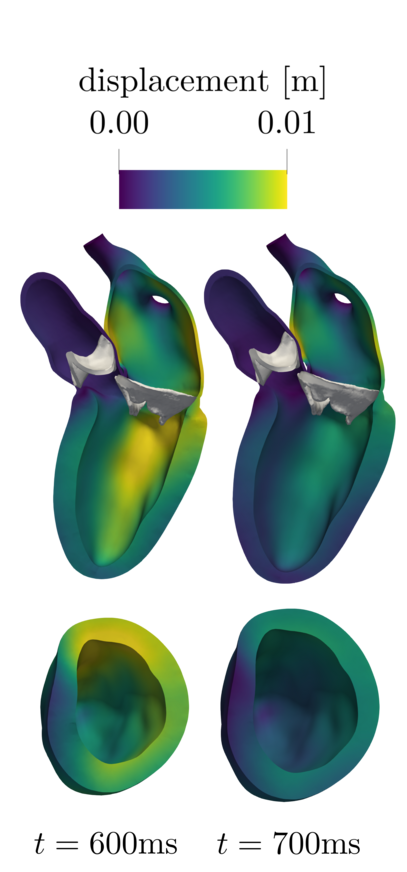}
        \caption{}
        \label{fig:filling-displacement}
    \end{subfigure}
    \begin{subfigure}{0.7596\figwidth}
        \centering
        \includegraphics[width=\textwidth]{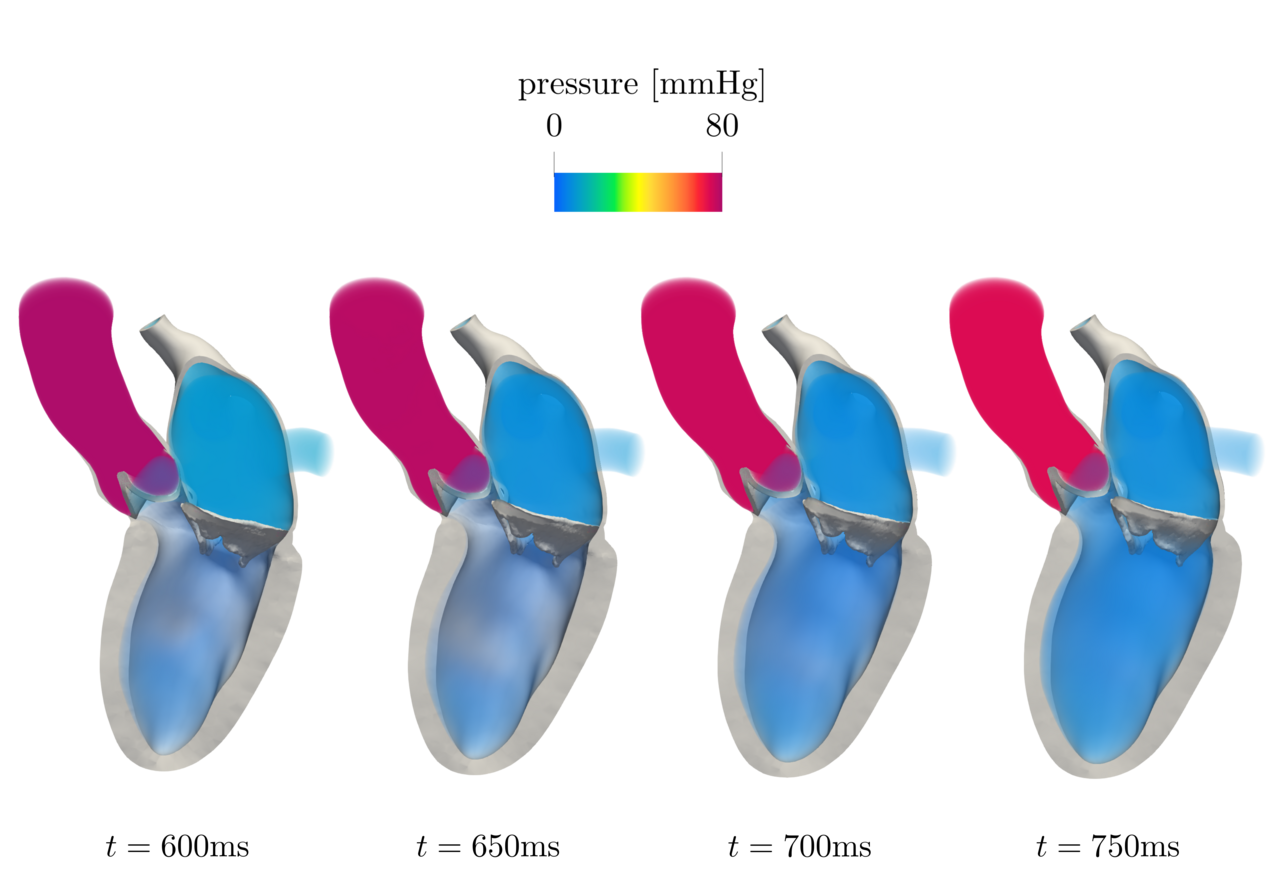}
        \caption{}
        \label{fig:filling-pressure}
    \end{subfigure}

    \caption{Solution snapshots during the filling phase, as seen through a long-axis and a short-axis section: (a) active tension $T_\text{act}$; (b) volume rendering of the fluid velocity magnitude $|\mathbf u|$; (c) displacement magnitude $|\mathbf d|$; (d) volume rendering of the fluid pressure $p$.}
    \label{fig:filling}
\end{figure}

\begin{figure}
    \centering
    \begin{subfigure}{0.25\textwidth}
        \includegraphics[width=\textwidth]{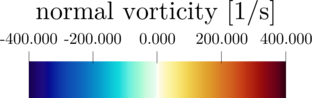}
    \end{subfigure}
    \vspace{0.5em}

    \begin{subfigure}{0.2\textwidth}
        \includegraphics[width=\textwidth]{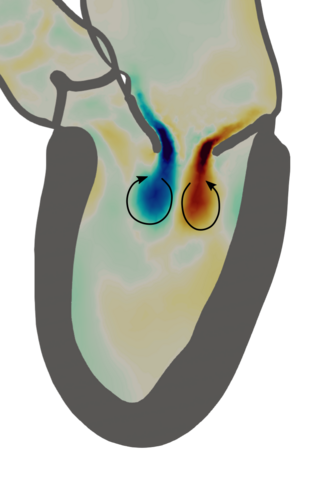}
        \caption{$t_\text{hb} = \SI{550}{\milli\second}$}
        \label{fig:solution-vorticity-a}
    \end{subfigure}
    \begin{subfigure}{0.2\textwidth}
        \includegraphics[width=\textwidth]{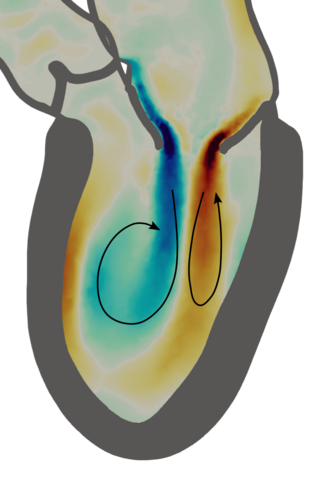}
        \caption{$t_\text{hb} = \SI{650}{\milli\second}$}
        \label{fig:solution-vorticity-b}
    \end{subfigure}
    \begin{subfigure}{0.2\textwidth}
        \includegraphics[width=\textwidth]{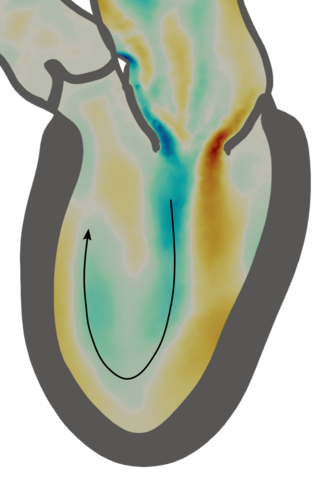}
        \caption{$t_\text{hb} = \SI{750}{\milli\second}$}
        \label{fig:solution-vorticity-c}
    \end{subfigure}

    \caption{Vorticity $\mathbf w\cdot\mathbf N = (\curl \mathbf u)\cdot\mathbf N$ of the velocity field, projected onto a slice of the domain (with normal $\mathbf N$ pointing outwards from the slice plane). The arrows indicate the direction of the rotating vortices.}
\end{figure}

Once the MV is open, blood flows from the atrium into the ventricle (\cref{fig:filling-velocity}). The flow is characterized by the formation of a jet through the MV, that is associated to a vortex ring (\cref{fig:solution-vorticity-a}). The vortex ultimately dissipates near to the ventricular free wall, while becoming larger near the septum (\cref{fig:solution-vorticity-b,fig:solution-vorticity-c}) resulting in a single vortex that rotates clockwise, if observed from a long-axis view with the septum on the left. This behavior is observed in medical images of healthy hearts \cite{di2018jet,pedrizzetti2014vortex}.
At the same time, the atrioventricular plane shifts towards the atrium, while the ventricular volume increases.
Since our model does not account for atrial contraction, the simulated diastolic phase lacks the atrial kick \cite{katz2010physiology}. This also determines higher-than-normal pressure and volume in the left atrium (with a peak pressure of \SI{19}{\mmhg}, against a normal value of about \SI{8}{\mmhg}) \cite{feher2017quantitative}. Filling continues until the ventricle starts contracting again, leading to the closure of the MV and the beginning of a new cardiac cycle.
The filling phase lasts for $T_\text{fil} = \SI{379}{\milli\second}$, and the whole diastolic phase lasts for $T_\text{dia} = \SI{474}{\milli\second}$, corresponding to \SI{59.3}{\percent} of the heartbeat.

\subsection{Conservation of blood volume}

\begin{figure}
    \centering
    \includegraphics[width=0.9\textwidth]{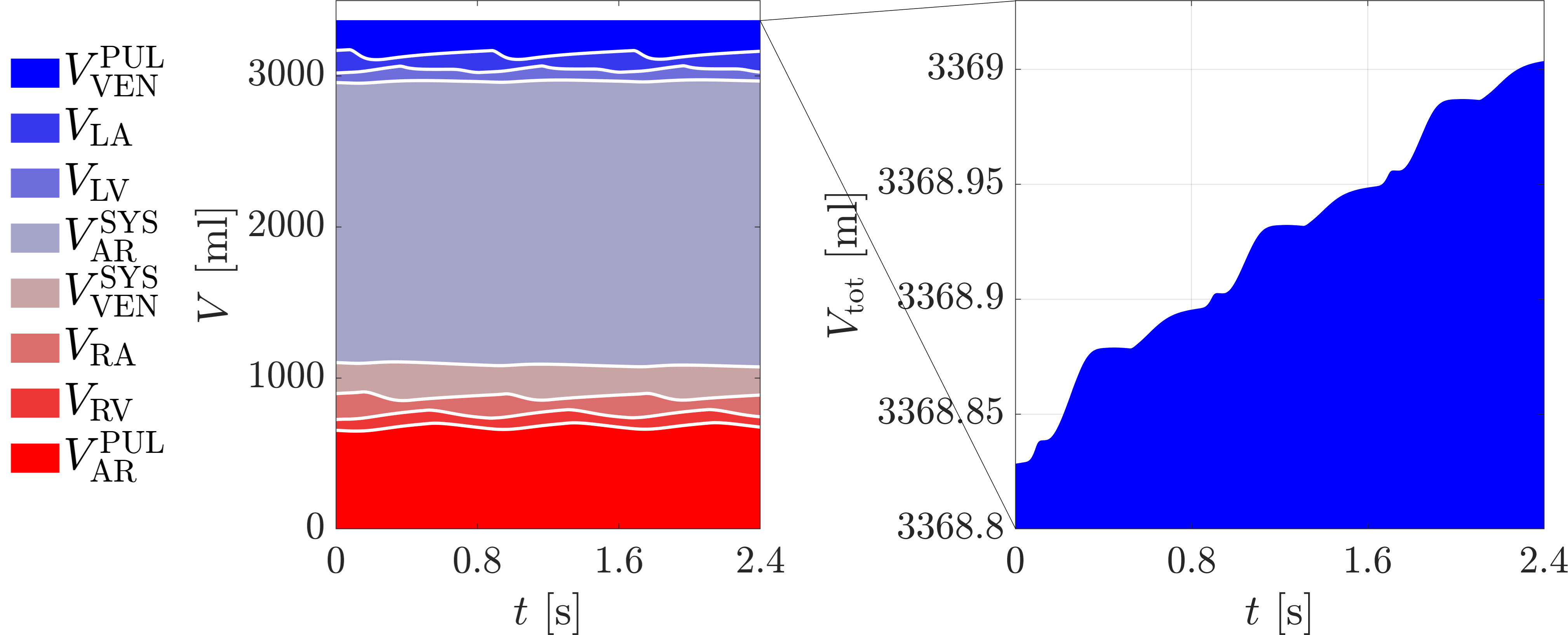}
    \caption{Blood volume over three heartbeats in each compartment of the heart and circulation. The zoom on the right shows how the total blood volume exhibits a very small variation in time: the mass gain over three heartbeats is only \SI{0.0052}{\percent} of the total.}
    \label{fig:blood-volume}
\end{figure}

The explicit treatment of the geometric FSI coupling condition and of the coupling between the Navier-Stokes equations and the circulation model might in principle lead to variations over time in the total blood volume. To assess whether this has an impact on simulation results, we compute the total blood volume over time as follows:
\begin{equation*}
    V_\text{tot}(t) = V\pul\ven(t) + V\la(t) + V\lv(t) + V\sys\ar(t) + V\sys\ven(t) + V_\text{RA}(t) + V_\text{RV}(t) + V\pul\ar(t)\;,
\end{equation*}
where $V\la(t)$ and $V\lv(t)$ are the volumes of $\Omega\la$ and $\Omega\lv\cup\Omega\ring$ at time $t$, respectively, and
\begin{align*}
    V\sys\ar(t) &= V\sys_{\text{AR},0} + V\ao(t) + C\sys\ar(p\sys\ar(t) - p_\text{EX}(t))\;, \\
    V\sys\ven(t) &= V\sys_{\text{VEN},0} + C\sys\ven(p\sys\ven(t) - p_\text{EX}(t))\;, \\
    V\pul\ar(t) &= V\pul_{\text{AR},0} + C\pul\ar(p\pul\ar(t) - p_\text{EX}(t))\;, \\
    V\pul\ven(t) &= V\pul_{\text{VEN},0} + C\pul\ven(p\pul\ven(t) - p_\text{EX}(t))\;.
\end{align*}
The latter equations are obtained by integrating the respective equations of the circulation model. We assume the zero-pressure volumes $V^\text{i}_\text{j,0}$ to be zero for simplicity, since they are constant in time and do not influence the assessment of blood volume conservation.

We report the blood volume over time for the different compartments, as well as the total blood volume, in \cref{fig:blood-volume}. From there, we can appreciate how the distribution of blood between the different compartments varies over time, but the total volume remains approximately constant. Indeed, the range of variation of the total volume over the three simulated heartbeats equals
\[
    \frac{\max_{t\in(0,T)} V_\text{tot}(t) - \min_{t\in(0,T)} V_\text{tot}(t)}{\max_{t\in(0,T)} V_\text{tot}(t)} = \SI{0.0052}{\percent}\;.
\]
We deem this very small variation over three heartbeats to be negligible and the result highly accurate. Therefore, the approximation introduced by the explicit discretization of geometric FSI coupling and circulation coupling does not introduce significant errors in terms of mass conservation.

\section{Limitations of the current study}
\label{sec:limitations}

In this section we present some limitations associated to the proposed computational framework, related to modeling choices and simplifications. To begin with, we modeled the LA as passive from both the electrical and mechanical viewpoints, for simplicity. On the contrary, in real hearts the left atrium is electrically active and contracts at the end of the diastolic phase, providing extra preload to the left ventricle along with an additional jet through the MV \cite{katz2010physiology,klabunde2011cardiovascular}. The inclusion of suitable atrial electromechanical models will be the subject of future studies.

The stimulation protocol used to trigger the activation of the myocardium is simplified if compared to the behavior of the cardiac conduction system \cite{costabal2016generating,del2022fast,landajuela2018numerical,katz2010physiology,romero2010effects,vergara2016coupled}. While this simplification is acceptable in physiological conditions \cite{regazzoni2022cardiac}, more detailed models might yield better descriptions of the ventricular activation pattern, which can become especially relevant if pathological scenarios are considered.

We employed the RIIS model for valve dynamics, prescribing the kinematics of the valve and neglecting any dynamic interaction with the blood. Moreover, we displace the valves along with the fluid domain displacement $\mathbf d\ale$, which has no physical meaning. Valve FSI plays a major role in the blood dynamics \cite{einstein2010fluid,feng2019analysis,fumagalli2021reduced,hsu2014fluid,luraghi2017evaluation,ma2013image,terahara2020heart}, as well as being responsible for several pathologies of clinical interest \cite{collia2019analysis,dabiri2019tricuspid}. While in this work we are mostly interested on the macroscopic effect of valves on the blood flow, the incorporation of suitable FSI models of the valves will be the subject of future studies.

Finally, in terms of numerical results, one major mismatch between our model and physiological data is the pressure difference between the ventricle and the aorta during the ejection phase. We obtain a maximum pressure difference of \SI{18}{\mmhg}. Such a high value, although not uncommon in the cardiac modeling literature \cite{this2020augmented,viola2020fluid}, is typically associated with a stenotic valve in clinical measurements \cite{jhun2021dynamics,wood1958aortic}. We believe this to be caused by minor inconsistencies in the geometrical model (based on a template obtained averaging over a sample of people, rather than on a specific individual) and by the simplicity of the valve model considered In any case, further investigations on this point are in order.

\section{Conclusions}
\label{sec:conclusions}

We introduced a novel fully integrated computational framework for the modeling and simulation of the human heart. The proposed model features three-dimensional and highly detailed descriptions of electrophysiology, active and passive mechanics and blood dynamics, as well as reduced models for valves and circulation. We included feedback effects among the different core models, resulting in a coupled multiphysics and multiscale problem.

Suitable numerical techniques were presented with the aim of solving efficiently the model equations. In particular, a staggered scheme with respect to time allowed for a flexible and efficient solution. Numerical methods were implemented in a high-performance computing framework within the \lifex library. We simulated a human left heart, and the results indicate that the proposed computational framework has the potential of describing accurately the physics of the heart, accounting for multiple feedback effects and overall capturing the interplay between the different processes driving the heartbeat. Our model is capable of representing all phases of the heartbeat and their durations, in accordance with clinical data. Moreover, through the coupling with a closed-loop circulation model, we were able to account for the interplay between the heart and the circulatory system, and guarantee the total blood volume conservation in time. The model correctly represents isovolumetric phases of the heartbeat, and reproduces effects associated to traveling pressure waves, thanks to the FSI coupling, suggesting that an EMF coupled modeling framework can overcome some of the limitations of electromechanics-driven CFD simulations. A more systematic comparison of the two modeling approaches will be the subject of future studies.

The agreement between numerical results and biomarkers from the medical literature suggests that the proposed model can be used for the investigation of the physiological cardiac function as well as for the simulation of pathological scenarios, and can serve as a stepping stone towards realistic and accurate integrated models of the whole heart.

\section{Acknowledgements}
This project has received funding from the European Research Council (ERC) under the European Union's Horizon 2020 research and innovation programme (grant agreement No 740132, iHEART - An Integrated Heart Model for the simulation of the cardiac function, P.I. Prof. A. Quarteroni).
\begin{center}
  \raisebox{-.5\height}{\includegraphics[width=.15\textwidth]{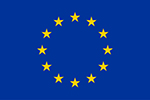}}
  \hspace{2cm}
  \raisebox{-.5\height}{\includegraphics[width=.15\textwidth]{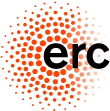}}
\end{center}

AZ received funding from the Italian Ministry of University and Research (MIUR) within the PRIN (Research projects of relevant national interest 2017 “Modeling the heart across the scales: from cardiac cells to the whole organ” Grant Registration number 2017AXL54F).

The authors acknowledge the CINECA award under the ISCRA B and ISCRA C initiatives, for the availability of high performance computing resources and support (ISCRA grants IsB22\_CoreMaS, P.I. Alfio Quarteroni, 2020-2022; IsC92\_HeartEMF, P.I. Michele Bucelli, 2021; IsC96\_EMFH, P.I. Michele Bucelli, 2022).

\appendix

\section{Solid mechanics constitutive laws}
\label{app:constitutive-laws}

The Guccione strain energy is computed from the displacement $\mathbf d$ as \cite{guccione1991finite}
\begin{equation*}
    \mathcal{W}_\mathrm{G}(\mathbf d) = \frac{c}{2}\left(\exp(Q) - 1\right) + \frac{\kappa}{2}(J - 1)\log(J)\,,
\end{equation*}
where
\begin{gather*}
    Q = \sum_{\mathrm{i, j} \in \{\mathrm{f}, \mathrm{s}, \mathrm{n}\}} a_\mathrm{i, j}\left(E \mathbf j\cdot \mathbf i\right)^2\;, \\
    E = \frac{1}{2}\left(F^T F - I\right)\;.
\end{gather*}

The strain energy function associated to the neo-Hookean model is given by \cite{ogden2013non}
\begin{equation*}
    \mathcal{W}_\mathrm{NH}(\mathbf d) = \frac{\mu}{2}\left(J^{-\frac{2}{3}}F:F - 3\right) + \frac{\kappa}{4}\left((J - 1)^2 + \log^2(J)\right)\;.
\end{equation*}

The values of the constant parameters appearing in both constitutive laws are reported in \cref{app:parameters}.

\section{Valve modeling}
\label{app:valves}

The heart model includes the MV in open configuration $\Gamma_\mathrm{MV}^\mathrm{open}$ and the AV in closed configuration $\Gamma_\mathrm{AV}^\mathrm{closed}$. The opposite configurations, $\Gamma_\mathrm{MV}^\mathrm{closed}$ and $\Gamma_\mathrm{AV}^\mathrm{open}$, are obtained in a preprocessing step as described in \cite{zingaro2022geometric}. The procedure also yields the displacement vectors between the two configurations, that is two fields $\mathbf d_\mathrm{k}: \Gamma_\mathrm{k}^\mathrm{closed} \to \mathbb{R}^3$ such that
\begin{equation*}
    \Gamma_\mathrm{k}^\mathrm{open} = \left\{ \mathbf x^\mathrm{open} \in \Omega\fluid\colon \mathbf x^\mathrm{open} = \mathbf x^\mathrm{closed} + \mathbf d_\mathrm{k}(\mathbf x^\mathrm{closed}),\;\mathbf x^\mathrm{closed} \in \Gamma_\mathrm{k}^\mathrm{closed} \right\}\;.
\end{equation*}
Then, at any time $t \in (0, T)$ the configuration $\Gamma_\mathrm{k}^t$ of the valve is given by
\begin{equation*}
    \Gamma_\mathrm{k}^t = \left\{ \mathbf x^\mathrm{t} = \mathbf x^\mathrm{closed} + c_\mathrm{k}(t)\,\mathbf d_\mathrm{k}(\mathbf x^\mathrm{closed}) + \mathbf d\ale(\mathbf x^\mathrm{closed}, t),\;\mathbf x^\mathrm{closed} \in \Gamma_\mathrm{k}^\mathrm{closed} \right\}\;,
\end{equation*}
wherein $c_\mathrm{k}(t)$ is a time-dependent opening coefficient that is equal to \num{0} when the valve is fully closed and to \num{1} when it is fully open. The functions $c_\mathrm{k}(t)$ are prescribed a priori, in the form
\begin{equation*}
    c_\mathrm{k}(t) = \begin{dcases}
        0 & \text{if } t \leq t^\mathrm{open}_\mathrm{k}\;, \\
        \frac{1}{2}\left(1 - \cos\left(\pi \frac{1 - \exp{-\chi_\mathrm{k}(t - t^\mathrm{open}_\mathrm{k})/\Delta t^\mathrm{open}_\mathrm{k}}}{1 - \exp(-\chi_\mathrm{k})}\right)\right) & \text{if } t^\mathrm{open}_\mathrm{k} < t \leq t^\mathrm{open}_\mathrm{k} + \Delta t^\mathrm{open}_\mathrm{k}\;, \\
        1 & \text{if } t^\mathrm{open}_\mathrm{k} + \Delta t^\mathrm{open}_\mathrm{k} < t < t^\mathrm{close}_\mathrm{k}\;, \\
        1 - \frac{1}{2}\left(1 - \cos\left(\pi \frac{1 - \exp{-\chi_\mathrm{k}(t - t^\mathrm{close}_\mathrm{k})/\Delta t^\mathrm{close}_\mathrm{k}}}{1 - \exp(-\chi_\mathrm{k})}\right)\right) & \text{if } t^\mathrm{close}_\mathrm{k} < t \leq t^\mathrm{close}_\mathrm{k} + \Delta t^\mathrm{close}_\mathrm{k}\;, \\
        0 & \text{if } t > t^\mathrm{close}_\mathrm{k} + \Delta t^\mathrm{close}_\mathrm{k}\;.
    \end{dcases}
\end{equation*}
The parameters $p_\mathrm{k}$, $\Delta t^\mathrm{open}_\mathrm{k}$ and $\Delta t^\mathrm{close}_\mathrm{k}$ are prescribed and they determine the opening speed of the valve. The opening time $t^\mathrm{open}_\mathrm{k}$ is the first instant at which the difference between upstream and downstream pressures is positive, while the closing time $t^\mathrm{close}_\mathrm{k}$ is the first instant at which the pressure difference is negative. Pressure difference is evaluated averaging over spherical control volumes as in \cite{zingaro2022geometric} (see \cref{fig:control-volumes}).

We remark that we displace the valves following the fluid domain. Moreover, we take $\mathbf u_{\Gamma_\mathrm{k}} = \mathbf 0$, i.e. we neglect the valve velocity due to the change of its configuration. This corresponds to a quasi-static approximation of the valve opening and closing \cite{fedele2017patient}.

\section{Blood circulation model}
\label{app:circulation}

The equations for the blood circulation model read: for $t \in (0, T)$,
\begin{equation}
    \begin{dcases}
        C\ar\sys\dv{p\ar\sys(t)}{t} = Q_\mathrm{AV}(t) - Q\ar\sys(t)\;, \\
        \frac{L\ar\sys}{R\ar\sys}\dv{Q\ar\sys(t)}{t} = -Q\ar\sys(t) - \frac{p\ven\sys(t) - p\ar\sys(t)}{R\ar\sys}\;, \\
        C\ven\sys\dv{p\ven\sys(t)}{t} = Q\ar\sys(t) - Q\ven\sys(t)\;, \\
        \frac{L\ven\sys}{R\ven\sys}\dv{Q\ven\sys(t)}{t} = -Q\ven\sys(t) - \frac{p_\mathrm{RA}(t) - p\ven\sys(t)}{R\ar\sys}\;, \\
        p_\mathrm{RA}(t) = p_\mathrm{EX}(t) + E_\mathrm{RA}(t)(V_\mathrm{RA}(t) - V_{0,\mathrm{RA}})\;, \\
        \dv{V_\mathrm{RA}(t)}{t} = Q\ven\sys(t) - Q_\mathrm{TV}(t)\;, \\
        Q_\mathrm{TV}(t) = \frac{p_\mathrm{RA}(t) - p_\mathrm{RV}(t)}{R_\mathrm{TV}(p_\mathrm{RA}(t), p_\mathrm{RV}(t))}\;, \\
        p_\mathrm{RV}(t) = p_\mathrm{EX}(t) + E_\mathrm{RV}(t)(V_\mathrm{RV}(t) - V_{0,\mathrm{RV}})\;, \\
        \dv{V_\mathrm{RV}(t)}{t} = Q_\mathrm{TV}(t) - Q_\mathrm{PV}(t)\;, \\
        Q_\mathrm{PV}(t) = \frac{p_\mathrm{RV}(t) - p\ar\pul(t)}{R_\mathrm{PV}(p_\mathrm{RV}(t), p\ar\pul(t))}\;, \\
        C\ar\pul\dv{p\ar\pul(t)}{t} = Q_\mathrm{PV}(t) - Q\ar\pul(t)\;, \\
        \frac{L\ar\pul}{R\ar\pul}\dv{Q\ar\pul(t)}{t} = -Q\ar\pul(t) - \frac{p\ven\pul(t) - p\ar\pul(t)}{R\ar\pul}\;, \\
        C\ven\pul\dv{p\ven\pul(t)}{t} = Q\ar\pul(t) - Q\ven\pul(t)\;, \\
        p_\text{LA}^\text{in}(t) = p\ven\pul(t) - R\ven\pul Q\ven\pul(t) - L\dv{Q\ven\pul(t)}{t}\;,
    \end{dcases}
\end{equation}
endowed with suitable initial conditions. In the above system, $p_\mathrm{EX}(t) = 0$ represents an external pressure, $E_\mathrm{RA}(t)$ and $E_\mathrm{RV}(t)$ are time-varying elastances modeling the contraction of the right atrium and ventricle, defined as
\begin{equation*}
E_\mathrm{i}(t) = E^\mathrm{B}_\mathrm{i} + E^\mathrm{B}_\mathrm{i} \phi_\mathrm{act}(t, t_\mathrm{c}^\mathrm{i}, T_\mathrm{c}^\mathrm{i}, T_\mathrm{r}^\mathrm{i})\;,
\end{equation*}
$\phi_\mathrm{act}$ is defined as in \cite{regazzoni2022cardiac}, and $V_{0,\mathrm{RA}}$ and $V_{0,\mathrm{RV}}$ are the resting volumes of left atrium and ventricle. Valvular resistances $R_\mathrm{TV}$ and $R_\mathrm{PV}$ are given by
\begin{equation*}
    R_\mathrm{k}(p_1, p_2) = \begin{dcases}
        R_{\max} & \text{if } p_1 < p_2\;, \\
        R_{\min} & \text{if } p_1 \geq p_2\;,
    \end{dcases} \qquad \mathrm{k} \in \{\mathrm{TV}, \mathrm{PV}\}\;.
\end{equation*}
The rest of the resistances, capacitances and inductances are parameters surrogating the properties of the circulation network.

\section{Model parameters}
\label{app:parameters}

\begin{table}
    \centering

    \begin{tabular}{l r S S}
        \toprule
        & \textbf{Parameter} & \multicolumn{2}{c}{\textbf{Value}} \\

        \midrule\multirow{2}{*}{Monodomain}
        & $\chi$ & 1400 & \si{\per\centi\metre} \\
        & $C_\text{m}$ & 1 & \si{\micro\farad\per\square\centi\metre} \\

        \midrule\multirow{3}{*}{Conductivities}
        & $\sigma_\mathrm{m}^\mathrm{l} / (\chi C_\text{m})$ & 2.00e-4 & \si{\square\metre\per\second} \\
        & $\sigma_\mathrm{m}^\mathrm{t} / (\chi C_\text{m})$ & 1.05e-4 & \si{\square\metre\per\second} \\
        & $\sigma_\mathrm{m}^\mathrm{n} / (\chi C_\text{m})$ & 0.55e-4 & \si{\square\metre\per\second} \\

        \midrule\multirow{3}{*}{Stimulus}
        & $A_\mathrm{app} / C_\text{m}$ & 25.71 & \si{\volt\per\second} \\
        & $\sigma_\mathrm{app}$ & 2.5e-3 & \si{\metre} \\
        & $T_\mathrm{app}$ & 3 & \si{\milli\second} \\
        \bottomrule
    \end{tabular}

    \caption{Parameters used in the electrophysiology model. Conductivities were tuned so as to obtain a conduction velocity of \SI{0.6}{\metre\per\second}, \SI{0.4}{\metre\per\second} and \SI{0.2}{\metre\per\second} along fibers, sheets and normal-to-fiber directions, respectively. The parameters used for the ionic model are those of the original paper \cite{ten2006alternans}.}
    \label{tab:params_ep}
\end{table}

\begin{table}
    \centering

    \begin{tabular}{r S S}
        \toprule
        \textbf{Parameter} & \multicolumn{2}{c}{\textbf{Value}} \\
        \midrule
        $\gamma$ & 30 & \\
        $k_d$ & 0.36 & \\
        $\alpha_{k_d}$ & -0.2083 & \\
        $K_\mathrm{off}$ & 8 & \si{\per\second} \\
        $K_\mathrm{basic}$ & 4 & \si{\per\second} \\
        $\mu_{fp}^0$ & 32.255 & \si{\per\second} \\
        $\mu_{fp}^1$ & 0.768 & \si{\per\second} \\
        $a_\mathrm{XB}$ & 8.9491e8 & \si{\pascal} \\
        $SL_0$ & 2.1 & \si{\micro\metre} \\
        \bottomrule
    \end{tabular}

    \caption{Parameters used in the force generation model. We only report those parameters whose values differ from the original setting described in \cite{regazzoni2020biophysically}.}
    \label{tab:params_activestress}
\end{table}

We report a list of the parameters used for the simulation described in \cref{sec:results}. \Cref{tab:params_ep} reports parameters for the electrophysiology model, \cref{tab:params_activestress} those for the force generation model, \cref{tab:params_mechanics} those for the solid mechanics model, \cref{tab:params_fluid} lists the parameters used in the fluid dynamics model and \cref{tab:params_circulation_a,tab:params_circulation_b} those of the circulation model. For the RDQ20-MF force generation model, we only report those parameters whose values differ from those presented in \cite{regazzoni2020biophysically}. We refer the interested reader to \cite{ten2006alternans} for details on the parameters of the TTP06 ionic model.

\begin{table}
    \centering

    \begin{tabular}{l r S S}
        \toprule
        & \textbf{Parameter} & \multicolumn{2}{c}{\textbf{Value}} \\

        \midrule
        \multirow{1}{3.75cm}{} & $\rho\solid$ & 1000 & \si{\kilo\gram\per\square\metre} \\

        \midrule\multirow{8}{3.75cm}{Guccione (atrium and ventricle)} &
        $c$ & 8.8e2 & \si{\pascal} \\
        & $a_\mathrm{ff}$ & 8 & \\
        & $a_\mathrm{ss}$ & 6 & \\
        & $a_\mathrm{nn}$ & 3 & \\
        & $a_\mathrm{fs}$ & 12 & \\
        & $a_\mathrm{fn}$ & 3 & \\
        & $a_\mathrm{sn}$ & 3 & \\
        & $\kappa$ & 5e4 & \si{\pascal} \\

        \midrule\multirow{2}{3.75cm}{Atrioventricular ring} &
        $\mu$ & 5e6 & \si{\pascal} \\
        & $\kappa$ & 1e6 & \si{\pascal} \\

        \midrule\multirow{2}{3.75cm}{Ascending aorta} &
        $\mu$ & 5.25e5 & \si{\pascal} \\
        & $\kappa$ & 1e6 & \si{\pascal} \\

        \midrule\multirow{4}{3.75cm}{Boundary conditions} &
        $K_\perp^\mathrm{epi}$ & 2e5 & \si{\pascal\per\metre} \\
        & $K_\parallel^\mathrm{epi}$ & 2e4 & \si{\pascal\per\metre} \\
        & $C_\perp^\mathrm{epi}$ & 2e4 & \si{\pascal\second\per\metre} \\
        & $C_\parallel^\mathrm{epi}$ & 2e3 & \si{\pascal\second\per\metre} \\

        \midrule\multirow{4}{3.75cm}{Initial conditions} &
        $p_0\la$ & 9.75 & \si{\mmhg} \\
        & $p_0\ring$ & 11.25 & \si{\mmhg} \\
        & $p_0\lv$ & 11.25 & \si{\mmhg} \\
        & $p_0\ao$ & 80 & \si{\mmhg} \\

        \midrule\multirow{2}{3.75cm}{Interface regularization}
        & $\psi_\text{th}\lv$ & 0.2 \\
        & $\psi_\text{th}\la$ & 0.1 \\

        \bottomrule
    \end{tabular}

    \caption{Parameters used in the solid mechanics model.}
    \label{tab:params_mechanics}
\end{table}

\begin{table}
    \centering

    \begin{tabular}{l r S S}
        \toprule
        & \textbf{Parameter} & \multicolumn{2}{c}{\textbf{Value}} \\

        \midrule\multirow{2}{*}{Navier-Stokes} &
        $\rho\fluid$ & 1060 & \si{\kilo\gram\per\cubic\metre} \\
        & $\mu\fluid$ & 3.5e-3 & \si{\pascal\second} \\

        \midrule\multirow{7}{*}{Valve modeling} &
        $R_\mathrm{MV}$, $R_\mathrm{AV}$ & 1e5 & \si{\kilo\gram\per\metre\per\second} \\
        & $\epsilon_\mathrm{MV}$, $\epsilon_\mathrm{AV}$ & 0.75e-3 & \si{\metre} \\
        & $\Delta t_\mathrm{MV}^\mathrm{open}$ & 10 & \si{\milli\second} \\
        & $\Delta t_\mathrm{AV}^\mathrm{open}$ & 10 & \si{\milli\second} \\
        & $\Delta t_\mathrm{MV}^\mathrm{close}$ & 30 & \si{\milli\second} \\
        & $\Delta t_\mathrm{AV}^\mathrm{close}$ & 80 & \si{\milli\second} \\
        & $\chi_\mathrm{MV}$, $\chi_\mathrm{AV}$ & -3 & \\
        \bottomrule
    \end{tabular}

    \caption{Parameters used in the fluid dynamics model. Valve half-thicknesses $\epsilon_\text{MV}$ and $\epsilon_\text{AV}$ were chosen to match literature data \cite{einstein2010fluid,sahasakul1988age}. Resistances $R_\mathrm{MV}$, $R_\mathrm{AV}$ were chosen to be sufficiently high to guarantee minimal spurious flow through valves without hindering the conditioning of the FSI system.}
    \label{tab:params_fluid}
\end{table}

\begin{table}
    \centering

    \begin{tabular}{l r S S}
        \toprule
        & \textbf{Parameter} & \multicolumn{2}{c}{\textbf{Value}} \\
        \midrule\multirow{6}{*}{Systemic arteries}
        & $R_\text{AR}^\text{SYS}$ & 0.45 & \si{\mmhg\second\per\milli\litre} \\
        & $C_\text{AR}^\text{SYS}$ & 2.19 & \si{\milli\litre\per\mmhg} \\
        & $L_\text{AR}^\text{SYS}$ & 2.7e-3 & \si{\mmhg\square\second\per\milli\litre} \\
        & $R_\text{upstream}^\text{SYS}$ & 0.07 & \si{\mmhg\second\per\milli\litre} \\
        & $p_\text{AR}^\text{SYS}(0)$ & 80 & \si{\mmhg} \\
        & $Q_\text{AR}^\text{SYS}(0)$ & 66.5775 & \si{\milli\litre\per\second} \\

        \midrule\multirow{5}{*}{Systemic veins}
        & $R_\text{VEN}^\text{SYS}$ & 0.26 & \si{\mmhg\second\per\milli\litre} \\
        & $C_\text{VEN}^\text{SYS}$ & 60 & \si{\milli\litre\per\mmhg} \\
        & $L_\text{VEN}^\text{SYS}$ & 5e-4 & \si{\mmhg\square\second\per\milli\litre} \\
        & $p_\text{VEN}^\text{SYS}(0)$ & 30.9029 & \si{\mmhg} \\
        & $Q_\text{VEN}^\text{SYS}(0)$ & 89.6295 & \si{\milli\litre\per\second} \\

        \midrule\multirow{6}{*}{Pulmonary arteries}
        & $R_\text{AR}^\text{PUL}$ & 0.05 & \si{\mmhg\second\per\milli\litre} \\
        & $C_\text{AR}^\text{PUL}$ & 10 & \si{\milli\litre\per\mmhg} \\
        & $L_\text{AR}^\text{PUL}$ & 5e-4 & \si{\mmhg\square\second\per\milli\litre} \\
        & $p_\text{AR}^\text{PUL}(0)$ & 20.0 & \si{\mmhg} \\
        & $Q_\text{AR}^\text{PUL}(0)$ & 69.3166 & \si{\milli\litre\per\second} \\

        \midrule\multirow{5}{*}{Pulmonary veins}
        & $R_\text{VEN}^\text{PUL}$ & 0.025 & \si{\mmhg\second\per\milli\litre} \\
        & $C_\text{VEN}^\text{PUL}$ & 38.4 & \si{\milli\litre\per\mmhg} \\
        & $L_\text{VEN}^\text{PUL}$ & 2.083e-4 & \si{\mmhg\square\second\per\milli\litre} \\
        & $p_\text{VEN}^\text{PUL}(0)$ & 17.0 & \si{\mmhg} \\
        & $Q_\text{VEN}^\text{PUL}(0)$ & 105.523 & \si{\milli\litre\per\second} \\

        \bottomrule
    \end{tabular}

    \caption{Parameters used in the circulation model: systemic and pulmonary circulation.}
    \label{tab:params_circulation_a}
\end{table}

\begin{table}
    \centering

    \begin{tabular}{l r S S}
        \toprule
        & \textbf{Parameter} & \multicolumn{2}{c}{\textbf{Value}} \\

        \midrule\multirow{2}{*}{Valves}
        & $R_\mathrm{min}$ & 7.5e-3 & \si{\mmhg\second\per\milli\litre} \\
        & $R_\mathrm{max}$ & 7.5e4  & \si{\mmhg\second\per\milli\litre} \\

        \midrule\multirow{8}{*}{Right atrium}
        & $E_\mathrm{A}$ & 0.06 & \si{\mmhg\per\milli\litre} \\
        & $E_\mathrm{B}$ & 0.07 & \si{\mmhg\per\milli\litre} \\
        & $t_\mathrm{C}$ & 0.8 & \\
        & $T_\mathrm{C}$ & 0.17 & \\
        & $T_\mathrm{R}$ & 0.17 & \\
        & $V_\mathrm{0,RA}$ & 4 & \si{\milli\litre} \\
        & $V_\mathrm{RA}(0)$ & 64.1702 & \si{\milli\litre} \\

        \midrule\multirow{8}{*}{Right ventricle}
        & $E_\mathrm{A}$ & 0.55 & \si{\mmhg\per\milli\litre} \\
        & $E_\mathrm{B}$ & 0.05 & \si{\mmhg\per\milli\litre} \\
        & $t_\mathrm{C}$ & 0.0 & \\
        & $T_\mathrm{C}$ & 0.34 & \\
        & $T_\mathrm{R}$ & 0.15 & \\
        & $V_\mathrm{0,RV}$ & 16 & \si{\milli\litre} \\
        & $V_\mathrm{RV}(0)$ & 148.9384 & \si{\milli\litre} \\

        \bottomrule
    \end{tabular}

    \caption{Parameters used in the circulation model: right heart.}
    \label{tab:params_circulation_b}
\end{table}

\clearpage
\bibliographystyle{abbrv}
\bibliography{bibliography}

\end{document}